\documentclass[12pt, a4paper]{article}
\usepackage[english]{babel}
\usepackage{amsmath}
\usepackage{amsthm}
\usepackage{amsfonts}
\usepackage{amssymb}
\usepackage{hyperref}
\usepackage{algorithm}
\usepackage{algpseudocode}
\usepackage{enumitem}
\usepackage{titlesec}
\usepackage{secdot}
\sectiondot{subsection}
\sectiondot{subsubsection}
\setlist[enumerate]{topsep=0pt, itemsep=.5ex, parsep=0ex}
\usepackage{tkz-tab}
\usepackage{fancyhdr}
\usepackage{color}
\usepackage[colorinlistoftodos]{todonotes}

\hypersetup{colorlinks=true, allcolors = blue}
\usepackage{tcolorbox}
\usepackage{subcaption}
\usepackage{titlesec}
\usepackage{relsize}
\usepackage{scalerel}
\titleformat*{\section}{\large\bfseries}

\tcbuselibrary{theorems}

\newtcbtheorem{mytheorem}{\emph{Định lý}}{colback = white!30, colframe = black!90, fonttitle = \bfseries}{th}
\usepackage{tabularx}
\usepackage[bitstream-charter]{mathdesign}

\usepackage[left=1.5cm,right=1.5cm,top=1.54cm,bottom=1.54cm,headheight=1cm]{geometry}
\usepackage{xcolor}
\hypersetup{
    colorlinks,
    linkcolor={red!70!black},
    citecolor={blue!70!black},
    urlcolor={blue!70!black}
}
\usepackage{mathtools}
\usepackage{multicol}

\usepackage{relsize}
\usepackage{fancyhdr}

\theoremstyle{remark}
\theoremstyle{definition}

\numberwithin{hq}{section}

\usepackage{pgf,tikz}
\usepackage{setspace}
\usepackage{mathrsfs}
\newcommand{\norm}[1]{\left\lVert#1\right\rVert}

\usepackage{makecell}
\usepackage{array, booktabs}
\usepackage[affil-it]{authblk}

\newcommand{\bu}{\pmb{u}}

\newcommand{\br}{\pmb{r}}

\newcommand{\bv}{\pmb{v}}
\newcommand{\bn}{\pmb{n}}

\newcommand{\bM}{\pmb{M}}

\newcommand{\bSig}{\pmb{\Sigma}}

\newcommand{\bK}{\pmb{K}}
\newcommand{\Ph}[1]{\textcolor{black}{#1}}
\newcommand{\Toan}[1]{\textcolor{black}{#1}} 

\begin{document}
\title{Parameter Estimation for the Reduced Fracture Model via a Direct Filter Method}
\author[1, *]{Phuoc Toan Huynh}
\author[1]{Feng Bao}
\author[2]{Thi-Thao-Phuong Hoang}
\affil[1]{Department of Mathematics, Florida State University, Florida 32304, USA}
\affil[2]{Department of Mathematics and Statistics, Auburn University, 221 Roosevelt Concourse, Auburn, 36849, Alabama, USA}
\affil[*]{Corresponding author. Email: th24ba@fsu.edu}

%

\maketitle
\abstract{In this work, we present a numerical method that provides accurate real-time detection for the widths of the fractures in a fractured porous medium based on observational data on porous medium fluid mass and velocity. To achieve this task, an inverse problem is formulated by first constructing a forward formulation based on the reduced fracture model of the diffusion equation. A parameter estimation problem is then performed online by utilizing a direct filter method. Numerical experiments are carried out to demonstrate the accuracy of our method in approximating the target parameters.
}

\section{Introduction}
Fluid flow and transport problems in fractured porous media have many important application in various fields such as subsurface hydrology, geophysics, and reservoir geomechanics. Therefore, it is necessary to have accurate and robust numerical simulations for such problems. However, such a task is often challenging due to the presence of the fractures. In particular, a fracture can represent either a fast pathway or a geological barrier, depending on whether its permeability is much higher or much lower than the surrounding rock matrix. In addition, the width of the fracture is much smaller than the size of the domain of calculation and any reasonable spatial mesh size. This issue requires refining the mesh locally around the fracture after imposing a global mesh on the entire domain, which is known to be computationally inefficient. One effective way to deal with this situation is to treat the fractures as interfaces, which avoids local refinement around the fractures. The original problem is then transformed into a new one where the interaction between the fractures and the surrounding rock matrix is taken into account (see~\cite{Alboin1999, Alboin2002, Amir2021, Angot2009, Fuma2011, Gander2021, Jaffre2005, SHLee2020, Morales2012} and the references therein). Models with such low-dimensional fractures are known as reduced fracture models or mixed-dimensional models. \Ph{Moreover, these models facilitate the use of different time scales in the fracture and in the rock matrix via a multi-scale approach, such as global-in-time domain decomposition methods}~(\cite{Hoang2013, Hoang2016, Hoang2017, Toan2023a, Toan2023b, Toan2024}). 

It is well-known that the characteristics of the fractures significantly influence fluid flow patterns and the selection of numerical algorithms. Thus, in this work, we consider the following inverse problem: Given the observations on the flow of the fluid in the fractured porous media, how can we infer the properties of the fractures, such as their widths? \Ph{This is a practically important problem} for engineers who seek to gain knowledge of the fractures when the observable data on the underground water is often limited. Such information will be invaluable when developing reliable models of naturally fractured reservoirs in order to optimize the primary and enhanced oil recovery methods. \Ph{The study of inversion schemes for state and parameter estimation in flow and transport in (fractured) porous media is challenging and has attracted great attention of researchers  (\cite{Ameur2002, Dai2004,LeGoc2010, Gwo2001, Krause2013,Mauldon1993, Neuman1979, Renshaw1996,  Wei2023}). However, limited work exists for the inverse problems associated with reduced fracture models. In \cite{Ameur2018}, an iterative algorithm was proposed for the steady-state problem to estimate the location and hydrogeological properties of a small number of lower-dimensional fractures in a porous medium using given distributed pressure or flow data.}

In this work, the fracture widths are considered as the parameters of interest, and the parameters detection problem we shall construct is based on a reduced fracture model of compressible fluid flow in which the fracture has larger permeability than the surrounding porous medium. One benefit of using the reduced fracture model is that one can include the fracture widths in the state estimation system, which facilitates \Ph{the} parameter estimation process. We aim to infer the widths of the fractures based on the real-time measurements of the porous media flow. To this end, the forward reduced fracture system is first implemented by using first-order backward Euler and mixed finite element methods~(\cite{Boffi2013,Brezzi1991,  Roberts1991}), which are mass conservative and can handle heterogeneous and anisotropic diffusion tensors effectively.  Next, to solve the inverse problem, we develop an online estimation method for the parameters representing the fractures widths. To achieve this goal, we utilize a direct filter method~(\cite{Bao2019a, Bao2021, Bao2023}) to dynamically estimate the unknown parameters as we receive the observational data in an online manner. The main idea of the direct filter method is to incorporate the physical model with the observations on the state process through the likelihood and use Bayesian inference to project the observational data onto the parameter space. Preliminary research has shown that the direct filter method can accurately estimate the parameters for high dimensional data assimilation~(\cite{Bao2022, Bao2019a}), and can be applied to solve practical problems~(\cite{Bao2021, Dyck2021, Bao2023}). \Ph{We remark that the model problem presented in this work is not complex, and our primary goal is to demonstrate the performance of the direct filter method for flows in fractured porous media. Thus, a multi-scale approach for the time discretization is not required here.}

The rest of this paper is organized as follows: In Section~\ref{sec2}, we introduce the model problem, formulate the corresponding reduced fracture model, and derive its fully discrete version. We shall discuss the optimal filtering for online parameter estimation in Section~\ref{sec3}. The direct filter method and its implementation are also introduced. In Section~\ref{sec4}, we discuss how to apply the direct filter method to estimate the target parameter in the reduced fracture model. Numerical experiments will also be carried out to illustrate the performance of the direct filter approach. The paper is then closed with a conclusion section.
\section{Reduced fracture model for the diffusion equations}\label{sec2}
\begin{figure}[h!]
\vspace{-0.4cm}
\centering
\includegraphics[scale=0.4]{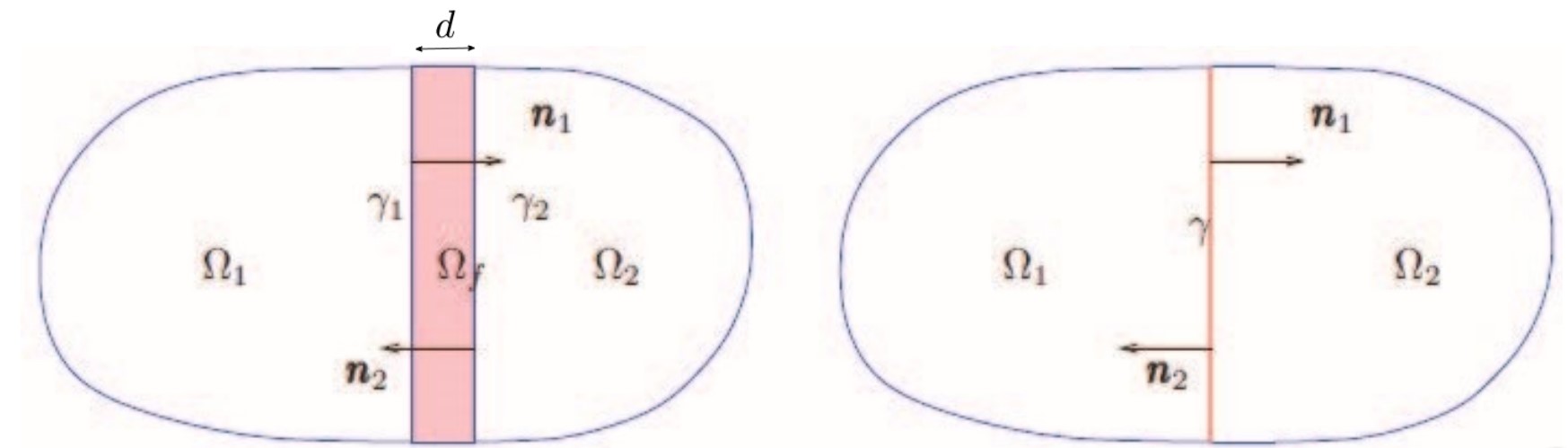}
\caption{The domain $\Omega$ with the fracture $\Omega_{f}$ (left) and the fracture-interface~$\gamma$~(right).}
\label{old_new_domain}  
\end{figure}

Let $\Omega$ be a bounded domain in $\mathbb{R}^{2}$ with Lipschitz boundary $\partial\Omega$, and $T>0$ be some fixed time.  Consider the flow problem of a single phase, compressible fluid written in mixed form as follows:
\begin{equation}
\label{original_problem}
\begin{array}{clll}
\phi\partial_t{p} + \text{div } \bu & = & q &\text{ in } \Omega \times (0, T), \\
\bu & = & -\bK \nabla{p} & \text{ in } \Omega \times (0, T), \\
p & = & 0 & \text{ on } \partial\Omega \times (0, T), \\
p(\cdot, 0) & = & p_0 & \text{ in } \Omega,
\end{array}
\end{equation}
where $p$ is the pressure, $\bu$ the velocity, $q$ the source term, $\phi$ the storage coefficient, and $\bK$ a symmetric, time-independent, hydraulic, conductivity tensor. Suppose that the fracture $\Omega_f$ is a subdomain of $\Omega$, whose thickness is $\Toan{d}$, that separates $\Omega$ into two connected subdomains: $\Omega \backslash \overline{\Omega}_f = \Omega_1 \cup \Omega_2, $ and $\Omega_1 \cap \Omega_2 = \emptyset.$ For simplicity, we assume further that $\Omega_f$ can be expressed as \vspace{-0.2cm}
$$
\Omega_f  = \left\{ \textbf{\textit{x}} \in \Omega: \textbf{\textit{x}} = \textbf{\textit{x}}_{\gamma} + s\bn \text{ where } \textbf{\textit{x}}_{\gamma} \in \gamma \text{ and } s \in \left(\Toan{-\dfrac{d}{2}, \dfrac{d}{2}}\right)\right\}, \vspace{-0.1cm}
$$
where $\gamma$ is the intersection between a line with $\Omega$ (see Figure~\ref{old_new_domain}). 

We denote by $\gamma_i$ the part of the boundary of $\Omega_i$ shared with the boundary of the fracture $\Omega_f$: $\, \gamma_i = (\partial\Omega_i \cap \partial\Omega_f) \cap \Omega$,  for $i =1,2$.  Let $\bn_i$ be the unit, outward pointing, normal vector field on $\partial\Omega_i$, where $\bn= \bn_1 = - \bn_2$. For $i=1,\; 2,\; f$, and for any scalar, vector, or tensor valued function $\varphi$ defined on $\Omega$, we denote by $\varphi_i$ the restriction of $\varphi$ to $\Omega_i$. The original problem \eqref{original_problem} can be rewritten as the following transmission problem: \vspace{-0.2cm}
\begin{equation}
\label{multidomain_problem}
\begin{array}{cllll}
\phi_i\partial_t{p_i} + \text{div }\bu_i & = & q_{i} & \text{ in } \Omega_i \times (0, T), & i = 1, 2, f, \\
\bu_i & = & - \bK_i \nabla{p}_i & \text{ in } \Omega_i \times (0, T), & i = 1, 2, f, \\
p_i & = & 0 & \text{ on } \left(\partial\Omega_i \cap \partial\Omega\right) \times (0, T), & i=1, 2, f, \\
p_i & = & p_f & \text{ on } \gamma_i \times (0, T), & i=1, 2, \\
\bu_i\cdot \bn_i & = & \bu_f\cdot\bn_i & \text{ on } \gamma_i \times (0, T), & i=1, 2, \\
p_i(\cdot, 0) & = & p_{0, i} & \text{ in } \Omega_i, & i=1, 2, f.
\end{array} \vspace{-0.2cm}
\end{equation}
The reduced fracture model that we consider in this paper was first proposed in~\cite{Alboin1999, Alboin2002} under the assumption that the fracture has larger permeability than that in the rock matrix. The model is obtained by averaging across the transversal cross sections of the two-dimensional fracture $\Omega_f$. We use the notation $\nabla_{\tau}$ and $\text{div}_{\tau}$ for the tangential gradient and tangential divergence, respectively. {We write $\phi_{\gamma}$ and $\bK_{\gamma}$ for $\Toan{d\phi_f}$ and $\bK_{f, \tau}$, respectively, where $\bK_{f, \tau}$ is the tangential component of $\bK_f$}. The reduced model consists of equations in the subdomains, \vspace{-0.2cm}
\begin{equation}
\label{reduced_subdomain}
\left.\begin{array}{rcll}
\phi_i\partial_t{p_i}+\text{div }\bu_i&=&q_{i} &\text{ in } \Omega_i\times (0, T), \\
\bu_i&=&-\bK_i\nabla{p_i} &\text{ in } \Omega_i\times (0, T), \\
p_i&=&0 &\text{ on } \left(\partial\Omega_i \cap \partial\Omega\right) \times (0, T), \\
p_i&=&p_{\gamma} &\text{ on } \gamma \times (0, T), \\
p_i(\cdot, 0)&=&p_{0, i} &\text{ in } \Omega_i, 
\end{array}\right.  \vspace{-0.2cm}
\end{equation}
for $i=1,2,$ and equations in the fracture-interface $\gamma$, \vspace{-0.2cm}
\begin{equation}
\label{reduced_fracture}
\begin{array}{rcll}
\phi_{\gamma}\partial_t{p_{\gamma}}+\text{div}_{\tau }\bu_{\gamma}&=& q_{\gamma} +\sum\limits^{2}_{i=1}\left( \bu_i \cdot \bn_i\right)_{\vert \gamma} & \text{ in } \gamma \times (0, T), \\
\bu_{\gamma} &=&\Toan{-\bK_{\gamma}d\nabla_{\tau}p_{\gamma}} & \text{ in } \gamma \times (0, T), \\
p_{\gamma}&=&0 &\text{ on } \partial\gamma \times (0, T), \\
p_{\gamma}(\cdot, 0)&=&p_{0, \gamma} & \text{ in } \gamma,
\end{array} \vspace{-0.2cm}
\end{equation}
where \Ph{$p_{\gamma}$, $\bu_{\gamma}$, and $q_{\gamma}$ are the reduced pressure, flux, and the source term, respectively, which are given by} 
\begin{align*}
\begin{array}{l}
\Ph{p_{\gamma} = \dfrac{1}{d}\mathlarger{\int}^{\frac{d}{2}}_{-\frac{d}{2}} p_{f}(x_{\gamma}+s\bn)ds}, \; \; \; \; \Ph{\bu_{\gamma} = \mathlarger{\int}^{\frac{d}{2}}_{-\frac{d}{2}}\bu_{f, \tau}(x_{\gamma}+s\bn)ds}, \; \; \; \; \Ph{q_{\gamma} = \mathlarger{\int}^{\frac{d}{2}}_{-\frac{d}{2}} q_{f}(x_{\gamma}+s\bn)ds}
\end{array}
\end{align*}
\Ph{where $\bu_{f, \tau}$ is the tangential component of $\bu_f$}. To write the weak formulation of \eqref{reduced_subdomain}-\eqref{reduced_fracture}, we use the convention that if $V$ is a space of functions, then $\pmb{V}$ is a space of vector functions having each component in $V$.  For arbitrary domain $\mathcal{O}$, we denote by $\left(\cdot, \cdot\right)_{\mathcal{O}}$ the inner product in $L^2\left(\mathcal{O}\right)$ or $\mathbf{\textbf{\textit{L}}^2\left(\mathcal{O}\right)}$. 
We next define the following Hilbert spaces: \vspace{-0.2cm}
\begin{align*}
\begin{array}{rl}
M &= \left\{v  = \left( v_1, v_2, v_{\gamma}\right) \in L^2\left(\Omega_1\right) \times L^2\left(\Omega_2\right) \times L^2\left(\gamma\right) \right\},  \vspace{0.2cm}\\
\Sigma &= \left\{\bv = \left(\bv_1, \bv_2, \bv_{\gamma}\right) \in \textbf{\textit{L}}^2\left(\Omega_1\right) \times \textbf{\textit{L}}^2\left(\Omega_2\right) \times \textbf{\textit{L}}^2\left(\gamma\right): \; \text{div} \; \bv_i \in L^2\left(\Omega_i\right), \; i=1,2, \right. \\
& \qquad \left.\text{and } \text{div}_{\tau} \; \bv_{\gamma} - \sum\limits^2_{i=1} \bv_i\cdot\bn_{i \vert \gamma} \in L^2(\gamma)\right\}.
\end{array}
\end{align*}
We define the bilinear forms $a(\cdot, \cdot; d)$, $b(\cdot, \cdot)$ and $c(\cdot, \cdot)$ on $\Sigma \times \Sigma$,  $\Sigma \times M$, and $M \times M$, respectively, and the linear form $L_q$ on $M$ by \vspace{-0.2cm}
\begin{equation}
\fontsize{10pt}{10pt}\selectfont
\label{BilinearForm}
\begin{array}{l}
\Toan{a\left(\bu, \bv; d \right)} = \sum\limits^2_{i=1} \left(\bK^{-1}_i\bu_i, \bv_i\right)_{\Omega_i} + \left(\left(\Toan{\bK_{\gamma}d}\right)^{-1}\bu_{\gamma}, \bv_{\gamma}\right)_{\gamma}, \hspace{0.1cm} b\left(\bu, \mu\right) = \sum\limits^2_{i=1}\left(\text{div}\; \bu_i, \mu_i\right)_{\Omega_i} + \left(\text{div}_{\tau} \; \bu_{\gamma} - \sum\limits^2_{i=1} \bu_i\cdot\bn_{i \vert \gamma}, \mu_{\gamma}\right)_{\gamma},\\
c_{\phi}(\eta, \mu) = \sum\limits^2_{i=1}\left(\phi_i\eta_i, \mu_i\right)_{\Omega_i} + \left(\phi_{\gamma}\eta_{\gamma}, \mu_{\gamma}\right)_{\gamma}, \quad L_{q}(\mu) = \sum\limits^2_{i=1}(q_i, \mu_i)_{\Omega_i}+(q_\gamma, \mu_\gamma)_{\gamma}. \vspace{-0.1cm}
\end{array}
\end{equation}
The weak form of \eqref{reduced_subdomain}-\eqref{reduced_fracture} can be written as follows: 

Find $p \in H^1(0, T; M)$ and $\bu \in L^2(0, T;\Sigma)$ such that \vspace{-0.1cm}
\begin{equation}
\label{weak_reduced}
\begin{array}{rcll}
\Toan{a\left(\bu, \bv; d \right)} - b\left(\bv, p \right) & = & 0 & \forall \bv \in \Sigma, \\
c_{\phi}\left(\partial_{t}p, \mu\right) + b\left(\bu, \mu\right) &=& L_q(\mu) & \forall \mu \in M,
\end{array} \vspace{-0.1cm}
\end{equation}
together with the initial conditions: \vspace{-0.1cm}
\begin{equation}
\label{initial_weak_reduced}
p_i(\cdot, 0)  =  p_{0, i}, \; \text{in} \; \Omega_i, \; i=1,2,  \quad \text{and} \quad p_{\gamma}(\cdot, 0)  =  p_{0, \gamma}, \; \text{in} \; \gamma.\vspace{-0.1cm}
\end{equation}
The well-posedness of problem \eqref{weak_reduced}-\eqref{initial_weak_reduced} was proved in~\cite{Hoang2016}.
To find the numerical solutions to \eqref{weak_reduced}-\eqref{initial_weak_reduced}, we discretize the problem in space using mixed finite element method (\cite{Boffi2013, Brezzi1991, Roberts1991}) and in time using backward Euler method. To this end, let $\mathcal{K}_{h, i}$ be a finite element partition of $\Omega_i \;(i=1, 2)$ into triangles. We denote by $\mathcal{G}_{h, i}$ the set of the edges of elements $\mathcal{K}_{h, i}$ lying on the interface~$\gamma$. Since $\mathcal{K}_1$ and $\mathcal{K}_2$ coincide on $\gamma$, the spaces $\mathcal{G}_{h, 1}$ and  $\mathcal{G}_{h, 2}$ are identical, thus we set $\mathcal{G}_h := \mathcal{G}_{h, 1} = \mathcal{G}_{h, 2}$. For $i=1, \, 2$, we consider the lowest order Raviart-Thomas mixed finite element spaces $M_{h, i} \times \Sigma_{h, i} \subset L^2(\Omega_i) \times H(\text{div}, \Omega_i)$:
\vspace{-0.1cm}
\begin{align*}
&M_{h, i}  =\left\{\mu_{h, i} \in L^2(\Omega_i): \; \mu_{h, i\vert K_i} = \text{const}, \; \; \forall K_i \in \mathcal{K}_{h, i} \right\}, \\
&\begin{array}{l}
\Sigma_{h, i} = \left\{\bv_{h, i} \in H(\text{div}, \Omega_i): \; \bv_{h, i\vert K_i} = (b_{K, i} + a_{K, i}x,  \; c_{K, i} + a_{K, i}y),\left(a_{K, i}, \; b_{K, i}, \; c_{K, i}\right) \in \mathbb{R}^3, \; \forall K_i \in \mathcal{K}_{h, i}\right\}.
\end{array}
\end{align*}
Similarly for the fracture, {let }$\Lambda_h \times \Sigma_{h, \gamma} \subset L^2(\gamma) \times H(\text{div}_{\tau}, \gamma)$ {be the lowest order Raviart-Thomas spaces in one dimension}:
\begin{align*}
& M_{h, \gamma} = \left\{ \mu_{h, \gamma} \in L^2(\gamma): \;  \lambda_{\vert E} = \text{const}, \; \forall E \in \mathcal{G}_h \right\}, \\
& \Sigma_{h, \gamma} = \left\{ \bv_{h, \gamma} \in H(\text{div}_{\tau}, \gamma): \; \bv_{h, \gamma\vert E} = az + b, \; (a, \;b) \in \mathbb{R}^2, \; \forall E \in \mathcal{G}_h \right\}.
\end{align*}

For the discretization in time, we consider a uniform partition of $\Ph{(0, \; T)}$ into $N$ sub-intervals $\Ph{(t^{n}, \; t^{n+1})}$ of length $\Delta{t}=t^{n+1}-t^{n}$, for $n=0, \hdots, N-1$. The time derivatives are approximated by the backward difference quotient
\begin{align*}
\bar{\partial}c^n = \dfrac{c^n-c^{n-1}}{\Delta{t}}, \; n= 1, \hdots, N,
\end{align*}
where the superscript $n$ indicates the evaluation of a function at the discrete time $t=t^n$. 

Finally, denote
\begin{align*}
\begin{array}{l}
\bM_h  = M_{h, 1} \times M_{h, 2} \times M_{h, \gamma}, \; \; \; \bSig_h = \Sigma_{h, 1} \times \Sigma_{h, 2} \times \Sigma_{h, \gamma},
\end{array}
\end{align*}
the fully-discrete version of \eqref{weak_reduced}-\eqref{initial_weak_reduced} reads as follows: 

For $n=1, \hdots, N$, find $(p^n_h, \bu^n_h) \in \bM_h \times \bSig_h$ satisfying
\begin{equation}
\label{weak_discrete_reduced}
\begin{array}{rcll}
\Toan{a\left(\bu^n_h, \bv_h; d\right)} - b\left(\bv_h, p^n_h \right) & = & 0 & \forall \bv_h \in \bSig_h, \vspace{0.2cm} \\
c_{\phi}\left(\bar{\partial}p^n_h, {\mu}_h\right) + b\left(\bu^n_h, {\mu}_h\right) &=& L_q^n({\mu}_h) & \forall \mu_h \in \bM_h,
\end{array} \vspace{-0.2cm}
\end{equation}
together with the initial conditions: \vspace{-0.1cm}
\begin{equation}
\label{initial_weak_discrete_reduced}
p^0_{h, i\vert K_i} := \dfrac{1}{\vert{K_i}\vert }\int_{K_i}p_{0, i}, \; \forall K_i \in \mathcal{K}_{h, i} \; i=1,2,  \quad \text{and} \quad p^0_{h, \gamma\vert E} = \dfrac{1}{\vert{E} \vert} \int_{E}p_{0, \gamma}, \; \forall E \in \mathcal{G}_h.\vspace{-0.1cm}
\end{equation}
\section{Data assimilation}\label{sec3}
In this section, we are interested in the following problems: given the observations on the solution of \eqref{weak_discrete_reduced}-\eqref{initial_weak_discrete_reduced},  what can we say about the characteristics of the existing fractures? In particular, we aim to estimate the widths of the fractures provided some partial information of the solution. To address such question, we need to incorporate some techniques from data assimilation.
\subsection{General framework}
We first briefly present the general framework of data assimilation. Suppose we have the following state-space model of a dynamical system
\begin{equation}
\label{StateEqs}
\begin{array}{l}
S_{n+1} = g(S_n) + \varepsilon_n, \; n=0, 1, \hdots,
\end{array}
\end{equation}
where $\Toan{g: \mathbb{R}^{l} \rightarrow \mathbb{R}^l}$ is a given mathematical model (linear of nonlinear), and the term $\varepsilon_n \in \Toan{\mathbb{R}^l}$ is the noise in the system which is usually assumed to be Gaussian. The process $\{S_n\}_n$ is referred to as the state process, which is often not directly observable, and the available data one received typically  consist of partial observations of the state $S_n$ which are given by
\begin{equation}
\label{ObservationEqs}
\begin{array}{l}
O_{n+1} = G(S_{n+1})+\zeta_{n+1}, \; n=0, 1, \hdots.
\end{array}
\end{equation}
The process $\{O_n\}_n \subset \mathbb{R}^k \; \Ph{(k \leq l)}$ is called the observation process, and $G: \Toan{\mathbb{R}^l} \rightarrow \mathbb{R}^k$ could be either a linear of nonlinear function, and $\zeta_{n+1}$ is a Gaussian noise. The goal of the data assimilation is to obtain the best estimate for the state $S_n$ given the observation process $\{O_i\}^n_{i=1}$. Mathematically, we need to find the optimal filer for $S_n$, denoted by $\bar{S}_n$ with 
\begin{equation}
\label{OptimalFilter}
\begin{array}{l}
\bar{S}_n: = \mathbb{E}\left[S_n \vert O_{1:n}\right].
\end{array}
\end{equation}
In this work, we use the Bayesian inference framework to obtain such estimation, which consists of two steps. The first step is the prediction step where the Chapman-Kolmogrov formula is applied to compute
\begin{equation}
\label{predict}
\begin{array}{l}
p\left(S_{n+1}\vert O_{1:n}\right) = \int p\left(S_{n+1}\vert S_n\right)p\left(S_n\vert O_{1:n}\right) dS_n,
\end{array}
\end{equation}
where we assume that $p\left(S_n\vert O_{1:n}\right)$ is known, and $p\left(S_{n+1}\vert S_n\right)$ is the transition kernel which can be computed by using \eqref{StateEqs}. 

The next step is to use the Bayes' formula to update the new posterior:
\begin{equation}
\label{update}
\begin{array}{l}
p\left(S_{n+1}\vert O_{1:{n+1}}\right) = \dfrac{p\left(O_{n+1}\vert S_{n+1}\right)p\left(S_{n+1}\vert O_{1:n}\right)}{p\left(O_{n+1}\vert O_{1:n}\right)}.
\end{array}
\end{equation}
As the given dynamical system is nonlinear in general, the terms occurring in \eqref{predict} and \eqref{update} may not have analytic expressions. The general idea to handle this difficulty is to construct an ensemble of particles, which represents a convex combination of Dirac measure, and iteratively update their locations and weights. As the number of particles increases, the conditional distribution is expected to converge to the exact distribution \Ph{(see, e.g.,~\cite{Doucet2001})}. This approach consists of the following steps: at time step $n$, assume we have a collection of particles $\{s^{(m)}_n\}^M_{m=1}$,
\begin{enumerate}
\item Approximate $p\left(S_n \vert O_{1:n}\right)$ by the empirical distribution $\tilde{p}\left(S_n \vert O_{1:n}\right)$ given by
\begin{equation}
\label{initial_empirical}
\begin{array}{l}
\tilde{p}\left(S_n \vert O_{1:n}\right) = \sum\limits^M_{m=1}w^{(m)}_n \delta_{s^{(m)}_n}\left(S_n\right)
\end{array}
\end{equation}
where $w^{(m)}_n$ stands for the weight of the particle $s^m_n$ which is given, and \Ph{$\delta_x$ is the Dirac delta function at $x$}.
\item For the prediction step, we first use the state equation \eqref{StateEqs} to generate $M$ empirical particles $\tilde{s}^{(m)}_{n+1}$ satisfying
\begin{align*}
\begin{array}{l}
\tilde{s}^{(m)}_{n+1} = g(s^{(m)}_n)+\varepsilon^{(m)}_n, \; m=1, \hdots, M.
\end{array}
\end{align*}
Those particles are used to compute the empirical distribution $\tilde{\pi}\left(S_{n+1}\vert O_{1:n}\right)$, which is an approximation for $p\left(S_{n+1}\vert O_{1:n}\right)$, through the Chapma-Kolmogrov formula:
\begin{equation}
\label{predict_empirical}
\begin{array}{l}
\tilde{\pi}\left(S_{n+1}\vert O_{1:n}\right) = \sum\limits^M_{m=1}w^{(m)}_n\delta_{\tilde{s}^{(m)}_{n+1}}\left(S_{n+1}\right).
\end{array}
\end{equation}
\item The empirical distribution of the update step is given by the following relation:
\begin{equation}
\label{update_empirical}
\begin{array}{l}
\tilde{\pi}\left(S_{n+1}\vert O_{1:{n+1}}\right) : = \sum\limits^M_{m=1}w^{(m)}_{n+1}\delta_{\tilde{s}^{(m)}_{n+1}}\left(S_{n+1}\right),
\end{array}
\end{equation}
where 
\begin{align*}
\begin{array}{l}
w^{(m)}_{n+1} = \dfrac{\tilde{w}^{(m)}_{n+1}}{\sum\limits^M_{m=1}\tilde{w}^{(m)}_{n+1}}, \vspace{0.1cm} \\
\tilde{w}^{(m)}_{n+1} = p\left(O_{n+1} \vert \tilde{s}^{(m)}_{n+1}\right)w^{(m)}_n.
\end{array}
\end{align*}
\item To avoid degeneracy, which means that a lot of the weights of the particles will be ignored, one needs a re-sampling step for the updated measure. This involves generating $M$ samples from the distribution \eqref{update_empirical} and assigning a weight of $\dfrac{1}{M}$ to each of the particle. The goal of this step is to eliminate particles with small weights.
\end{enumerate}

Since estimating the parameter $\Toan{d}$ is our main focus and the state estimates are unnecessary, we adopt the Direct Particle Filter method developed in~\cite{Bao2019a}. Such method will help us avoid unwanted distraction in the state estimation process and directly provide the approximations of the target parameter. As a consequence, the dimension of the problem is the same as the dimension of the parameter of interest, which would address the "curse of dimensionality" when solving nonlinear filtering problems in online parameter estimation (\cite{Bao2019a}). Therefore, this method could solve the parameter estimation problem efficiently, especially for our case when the dimension of the state model is very high and the parameter $\Toan{d}$ is a low dimensional vector.
\subsection{Direct filter for parameter estimation}
Similar to \eqref{StateEqs}-\eqref{ObservationEqs}, let $X_n \in \mathbb{R}^l$ describe the state of some physical model, and $Y_n \in \mathbb{R}^k$ is a noisy observation of $X_n$ with a noise perturbation $\xi_n$. We also denote by $\theta \in \mathbb{R}^p \; (p \ll l)$ the vector of target parameters, and assume that we have the following system for the parameter estimation problem:
\begin{align}
X_{n+1} = h(X_n, \;\theta) + w_n,\label{ParaEts_Eqs} \vspace{0.1cm} \\
Y_{n+1} = HX_{n+1} + \xi_{n+1},\label{ObservationalPara_Eqs}
\end{align}
where $h: \Toan{\mathbb{R}^l} \times \mathbb{R}^p \rightarrow \Toan{\mathbb{R}^l}$ is a function (linear or nonlinear) representing the considered physics model, $H: \Toan{\mathbb{R}^l} \rightarrow \mathbb{R}^k$ is a linear matrix, and $(w_n, \xi_n)$ is a pair of independent Gaussian noises.  The parameter estimation problem we are interested in is to estimate $\theta$ in \eqref{ParaEts_Eqs} by using the observational data $Y$ provided in \eqref{ObservationalPara_Eqs}.

In the direct filter method, instead of treating $\theta$ as a vector containing deterministic constant entries, we consider $\theta$ as a stochastic process to be estimated with respect to time. To this end, we replace $\theta $ in $\eqref{ParaEts_Eqs}$ with $\theta_n$, and rewrite \eqref{ParaEts_Eqs}-\eqref{ObservationalPara_Eqs} in the following form:
\begin{align}
&\theta_{n+1} = \theta_{n}+\epsilon_n, \vspace{0.1cm} \\ 
&Y_{n+1} = H\left(h(X_n, \; \theta_{n+1}) + w_n\right) + \xi_{n+1}, 
\end{align}
where $\epsilon_n$ is an artificial dynamic noise. We define $\zeta_{n+1} = Hw_n+\xi_{n+1}$, which is a multivariate Gaussian variable, and we can obtain the following dynamics
\begin{align}
&\theta_{n+1} = \theta_{n}+\epsilon_n,\label{ParaStos_Eqs} \vspace{0.1cm} \\ 
&Y_{n+1} = Hh(X_n, \;\theta_{n+1})+ \zeta_{n+1}. \label{NewObser_Eqs}
\end{align}
We aim to find the best estimate $\mathbb{E}\left[\theta_n \vert Y_{1:n}\right]$ for $\theta_n$ at any time instant $n$ given the observation process $\{Y_i\}^n_{i=1}
$. Instead of generating a long-term simulation trajectory for the process of $X_n$  (e.g.~\cite{Kantas2015}), we follow the idea in~\cite{Bao2019a} and use the fact that the observation data $Y_n$ provides direct observations on $X_n$ and introduce the following approximation scheme
\begin{align*}
H^{-1}Y_n \approx H^{-1}(Y_n - \xi_n) = X_n,
\end{align*}
which leads to the approximation dynamical system
\begin{align}
&\theta_{n+1} = \theta_{n}+\epsilon_n,\label{ParaStos_ApproEqs} \vspace{0.1cm} \\ 
&\tilde{Y}_{n+1} = Hh(H^{-1}\tilde{Y}_n, \; \theta_{n+1})+ \zeta_{n+1}. \label{NewObser_ApproEqs}
\end{align}
We conclude this section by introducing the implementation of the direct filter \Ph{method}. More detailed discussions can be found in~\cite{Bao2019a}. This method is initialized by first generating a collection of $M$ particles $\left\{\Ph{\theta^{(m)}_n}\right\}^M_{m=1}$ that describes $p\left(\theta_n \vert Y_{1:n}\right)$ at time instant $n$ and consists of the following three steps:
\begin{enumerate}
\item \textbf{The prediction step.} Generates a prior estimate $\tilde{\theta}_{n+1} = \left\{\tilde{\theta}^{(m)}_{n+1}\right\}^M_{m=1}$ for the target parameter $\theta_{n+1}$ using \eqref{ParaStos_ApproEqs} by adding $\left\{\epsilon^{(m)}_n\right\}^M_{m=1}$ to $\left\{\theta^{(m)}_n\right\}^M_{m=1}$:
\begin{align*}
\begin{array}{l}
\tilde{\theta}^{(m)}_{n+1} = \theta^{(m)}_{n} + \epsilon^{(m)}_n, \; m =1, 2, \hdots, M.
\end{array}
\end{align*}
The prediction step provides an empirical distribution $\tilde{\pi}\left(\theta_{n+1}\vert Y_{1:n}\right)$ for the prior $p\left(\theta_{n+1} \vert Y_{1:n}\right)$ which is computed by:
\begin{equation}
\label{DirectPredict}
\begin{array}{l}
\tilde{\pi}\left(\theta_{n+1} \vert Y_{1:n}\right) = \dfrac{1}{M}\sum\limits^M_{m=1} \delta_{\tilde{\theta}^{(m)}_{n+1}}(\theta_{n+1}).
\end{array}
\end{equation}
\item \textbf{The update step.} The update step incorporates the observational data and derives a weighted posterior distribution based on the prior $\tilde{\pi}\left(\theta_{n+1} \vert Y_{1:n}\right)$ given in \eqref{DirectPredict}:
\begin{equation}
\label{DirectUpdate}
\begin{array}{l}
\tilde{\pi}\left(\theta_{n+1} \vert Y_{1:n+1}\right) = \dfrac{1}{C}\sum\limits^M_{m=1}p\left(Y_{n+1}\vert \tilde{\theta}^{(m)}_{n+1}\right)\delta_{\tilde{\theta}^{(m)}_{n+1}}(\theta_{n+1}) = \sum\limits^M_{m=1} w^{(m)}_{n+1}\delta_{\tilde{\theta}^{(m)}_{n+1}}(\theta_{n+1}),
\end{array}
\end{equation}
where $w^{(m)}_{n+1} = p\left(Y_{n+1}\vert \tilde{\theta}^{(m)}_{n+1}\right)/C$ is the weight for the particle $\tilde{\theta}^{(m)}_{n+1}$ with a normalization factor $C$ and the likelihood function in \eqref{DirectPredict} is given by
\begin{align}
\begin{array}{l}
p\left(Y_{n+1}\vert \tilde{\theta}^{(m)}_{n+1}\right) = \text{exp}\left(-\dfrac{1}{2}\norm{Hh\left(H^{-1}Y_n, \tilde{\theta}^{(m)}_{n+1}\right)- Y_{n+1}}^2_R\right),
\end{array}
\end{align}
where $\norm{\alpha}_{R}:= \alpha R^{-1}\alpha$ with $R$ standing for the invariance variance of the observational noise $\Toan{\zeta_{n+1}}$.
\item \textbf{Resampling step.} The purpose of the resampling step is to generate a set of equally weighted samples to avoid the degeneracy issue (see~\cite{Bao2019a, Bao2019b, Bao2014, Bao2020, Bao2018, Bao2017}). In this work, we simply use the importance sampling method~(\cite{Doucet2001b, Morzfeld2018}) to generate samples, denoted by $\left\{\theta^{(m)}_{n+1}\right\}^M_{m=1}$, from the weighted importance distribution $\tilde{\pi}\left(\theta_{n+1}\vert Y_{1:{n+1}}\right)$. \Ph{More specifically, we use the distribution $\tilde{\pi}\left(\theta_{n+1} \vert Y_{1:n+1}\right)$ to sample with replacement $M$ particles $\left\{\tilde{\theta}^{(m)}_{n+1}\right\}^M_{m=1}$ according to the weights $\left\{w^{(m)}_{n+1}\right\}^M_{m=1}$ to obtain $\left\{\theta^{(m)}_{n+1}\right\}^M_{m=1}$. This way, we can replace particles with small weights by particles with large weights. We then obtain an ``unweighted" empirical distribution as follows}
\begin{align*}
\begin{array}{l}
\Ph{\tilde{p}\left(\theta_{n+1} \vert Y_{1:n+1}\right) = \dfrac{1}{M}\sum\limits^M_{m=1}\delta_{{\theta}^{(m)}_{n+1}}(\theta_{n+1}),}
\end{array}
\end{align*}
\Ph{and use $\tilde{p}\left(\theta_{n+1} \vert Y_{1:n+1}\right)$ as the initial empirical distribution for the next time instant $n+1$.}
\end{enumerate}

With the above procedure, the estimate we obtain for the target parameter at the time instant $n+1$ is given by
\begin{align*}
\tilde{\theta} = \dfrac{1}{\Toan{n+1-j}}\sum\limits^{n+1}_{\Toan{i=j}} \tilde{\mathbb{E}}\left[\theta_i \vert Y_{1:i}\right],
\end{align*}
where $\Toan{j}$ is a number of burn-in steps to reduce the influence of large noises at some time instants, and $\tilde{\mathbb{E}}\left[\theta_i \vert Y_{1:i}\right]$ is the mean estimate of the empirical distribution $\Toan{\tilde{p}\left(\theta_i\vert Y_{1:i}\right)}$.
\section{Numerical results}\label{sec4}
In this section, we shall explain how to apply the direct filter method to estimate the width of the fracture based on the fully discrete version of the reduced fracture model, i.e., the system of equations \eqref{weak_discrete_reduced}-\eqref{initial_weak_discrete_reduced}. We then carry out some numerical experiments with different fracture configurations to demonstrate the performance of the method.
\subsection{Problem setup}
To incorporate the direct filter algorithm, we first rewrite the system \eqref{weak_discrete_reduced}-\eqref{initial_weak_discrete_reduced} in a more compact form. Let $\bar{\pmb{M}}_h$ and $\bar{\bSig}_h$ be the finite collections of basis functions in $\pmb{M}_h$ and $\bSig_h$, respectively. We introduce the following matrices resulted in the bi-linear forms given in \eqref{BilinearForm}:
\begin{align}
\begin{array}{l}
\pmb{A}_h(d) = \left(a(\br_h, \bv_h; d)\right)_{\br_h, \bv_h \in \bar{\bSig}_h}, \; \pmb{B}_h = \left(-b(\bv_h, {\eta}_h)\right)_{\bv_h\in \bar{\bSig}_h, \eta_h \in \bar{\pmb{M}}_h}, \; \pmb{C}_{h, \phi} = \left(-c_{\phi}({\eta}_h, {\lambda}_h)\right)_{{\eta}_h, {\lambda}_h \in \bar{\pmb{M}}_h}. 
\end{array}
\end{align}
By replacing $v_h$ and $\mu_h$ in \eqref{weak_discrete_reduced} by the basis functions and expressing $\bu^n_h$ and $p^n_h$ in terms of those basis functions, we obtain the following matrix form for \eqref{weak_discrete_reduced}
\begin{equation}
\label{MatForm_v1}
\begin{array}{l}
\pmb{\Lambda}_h(d)X_n = F(X_{n-1}, \; \pmb{G}_h),
\end{array}
\end{equation}
where
\begin{align}
\begin{array}{l}
\pmb{\Lambda}_h(d) =\begin{bmatrix}
\pmb{A}_h(d) & \pmb{B}_h \vspace{0.1cm} \\
\Delta{t}\left(\pmb{B}_h\right)^T  &\pmb{C}_{h, \phi},
\end{bmatrix}
\end{array}
\end{align}
with $\Delta{t}$ being the chosen timestep size, $\pmb{G}_h$ is a vector containing the boundary conditions, and $X_n = \left(\bu^{n}_h, p^n_h\right)$ . Equivalently, we have the following recursive relation
\begin{equation}
\label{StateEqs_Discrete}
\begin{array}{l}
X_n = \pmb{\Phi}(X_{n-1}, \; d), \; n=1, 2, \hdots,
\end{array}
\end{equation} 
where 
\begin{align*}
\pmb{\Phi}(X_{n-1}, d) = \left(\pmb{\Lambda}(d)\right)^{-1}F(X_{n-1}, \;\pmb{G}_h).
\end{align*}
We can see that by using the reduced fracture model, the parameter of interest $d$ is included in the state equation \eqref{StateEqs_Discrete}. Thus, the direct filter method can be applied to approximate $d$.  For all test cases shown in this section, we fix $\Delta{t} = 0.1$ \Ph{and the spatial meshsize $h = 1/50$}, therefore, we have a total collection of $K = T/\Delta{t}$ such state vectors, where $T$ is the terminal time.

To formulate the data assimilation problem for parameter estimation, we assume that we receive datasets
\vspace{-0.1cm}
\begin{align*}
\begin{array}{l}
Y_{n+1} = HX_{n+1} + \xi_{n+1},
\end{array}
\vspace{-0.1cm}
\end{align*}
as the observational process with linear dependence on $X_{n+1}$ and $\xi_{n+1}$ is a noise following a Gaussian distribution. As we aim to learn the information about the fracture, we assume that the received data is directly related to the solution on the fracture. In particular, we have
\begin{align*}
\begin{array}{l}
HX_{n+1}= X_{\gamma, n+1},
\end{array}
\vspace{-0.1cm}
\end{align*}
where $X_{\gamma, n+1} = \left(\bu^{n+1}_{h, \gamma}, p^{n+1}_{h, \gamma}\right)$ is the solution on the fracture at time instant $n+1$. As the width of the fracture in general is significantly small compared to the overall size of the domain of calculation, instead of approximating $d$, we denote $\Toan{\theta := 1/d}$ to be the parameter of interest, and we add noise to $\theta$ to make it a stochastic process. The introduction of additional noise to the parameter of interest will transform it into an ensemble. As such, data assimilation techniques such as particle methods can be applied to facilitate the state estimation of those parameters. All together, we obtain the following optimal filtering problem for the parameter estimation task:
\begin{align}
& \theta_{n+1} = \theta_{n} + \epsilon_n, \label{ParaSetup_Egs} \vspace{0.1cm} \\
& Y_{n+1} = H\pmb{\Phi}\left(H^{-1}(Y_n), \;1/\theta_{n+1}\right) + \xi_{n+1}. \label{ObserSetup_Egs}
\end{align}

Note that \eqref{ParaSetup_Egs}-\eqref{ObserSetup_Egs} is derived from the system \eqref{weak_discrete_reduced}-\eqref{initial_weak_discrete_reduced} for the case with a single fracture. Similar reduced fracture models can be extended readily to the other test cases presented in this section, and the corresponding parameter estimation system can then be formulated. In what follows, we consider three different test cases with different settings of the fractures to study the applicability of the direct filter method to estimate the parameter of interest.
\subsection{Test case 1: One Single Fracture} 
\begin{figure}[h]
\centering
\vspace{-0.3cm}
\begin{minipage}{.5\textwidth}
  \centering
  \includegraphics[width=0.85\linewidth]{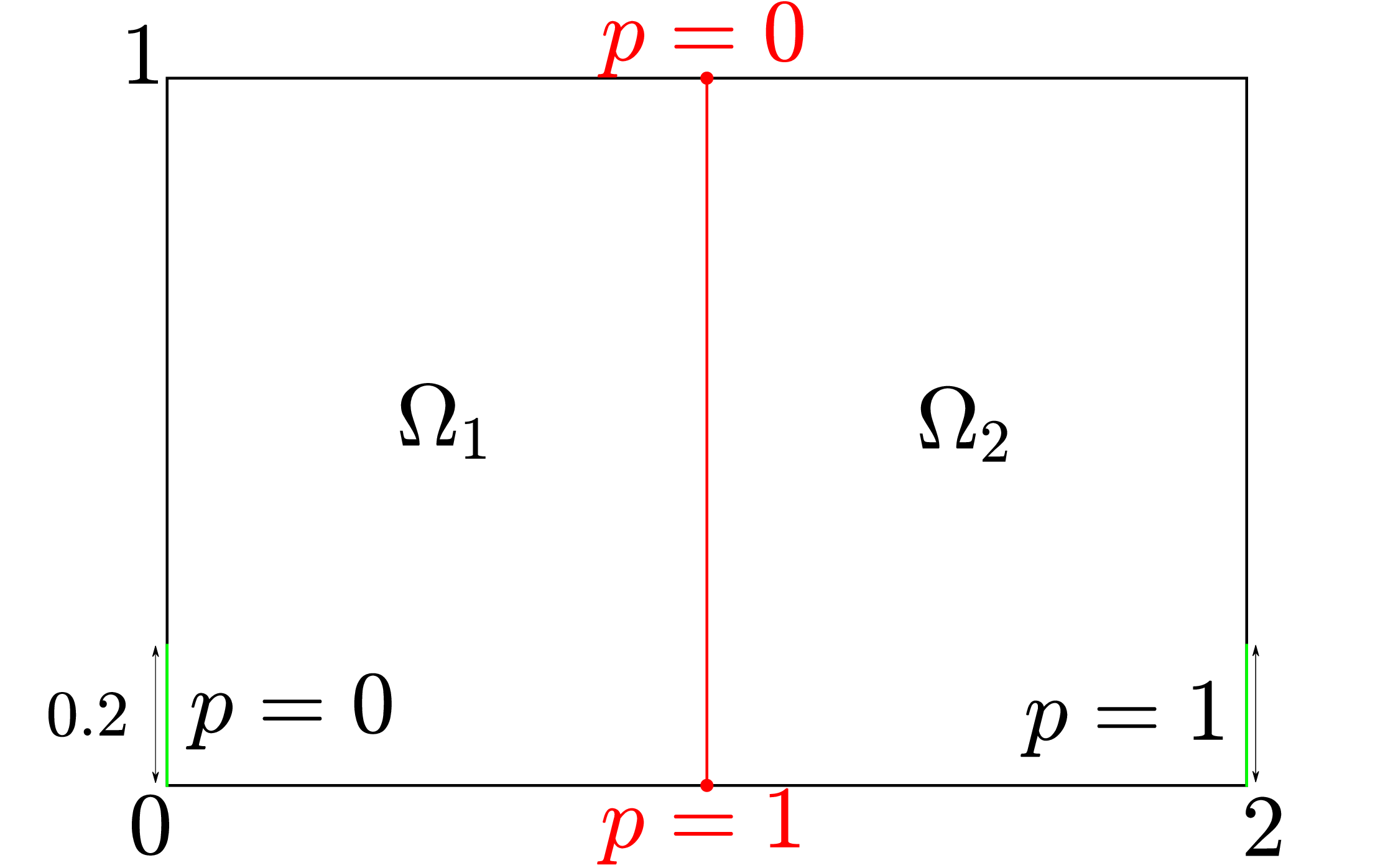}
  \label{fig:test1}
\end{minipage}%
\begin{minipage}{.5\textwidth}
  \centering
  \includegraphics[width=0.85\linewidth]{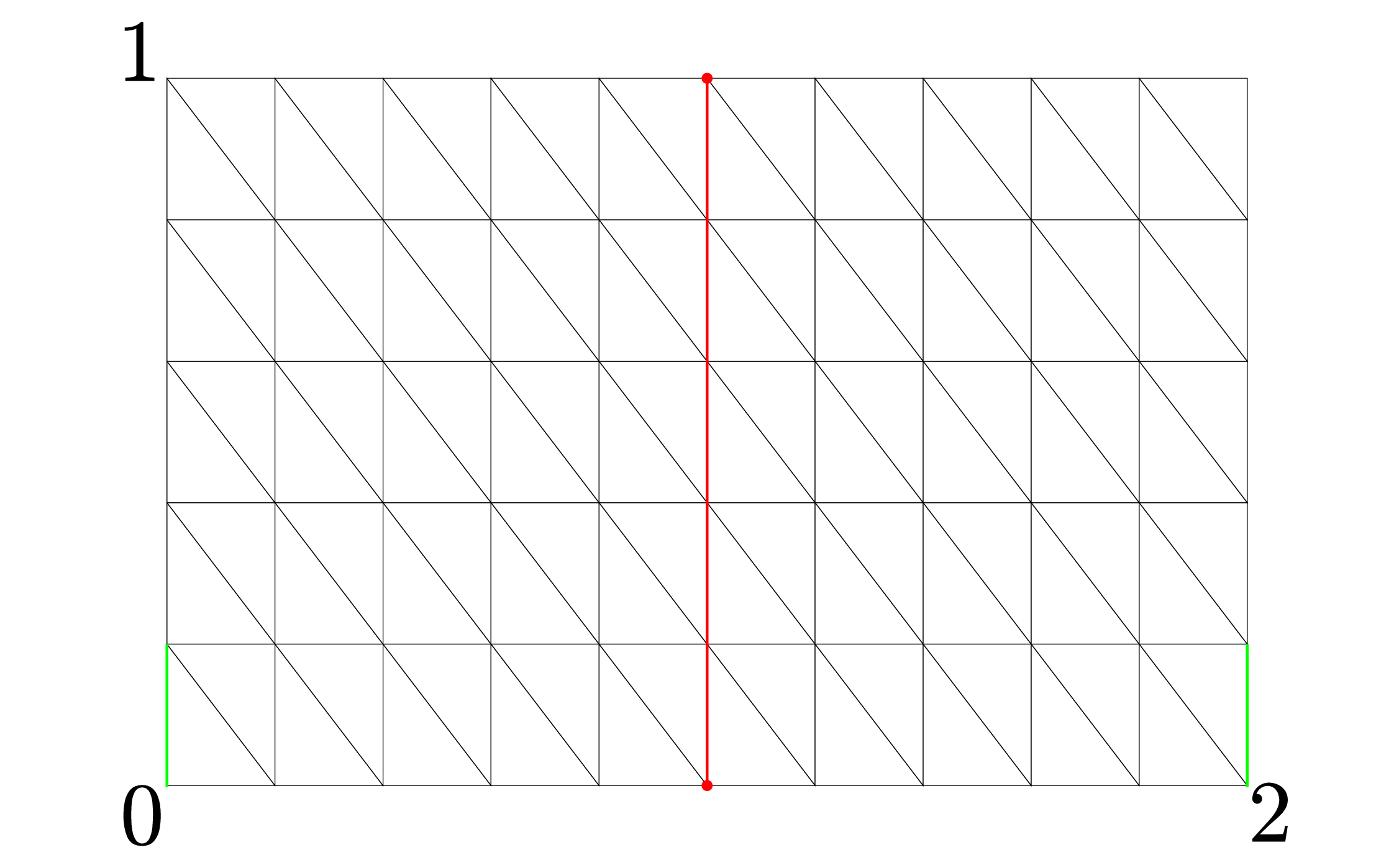}
  \label{fig:test2}
\end{minipage} 
\caption{[Test case 1] (Left) Geometry and boundary conditions of the test case.  (Right) Example of an uniform triangular mesh for spatial discretization.}
\label{Non_Immersed} 
\vspace{-0.4cm}
\end{figure}

\begin{figure}[h!]
\centering
\begin{minipage}{.47\textwidth}
  \centering
  \includegraphics[width=1\linewidth]{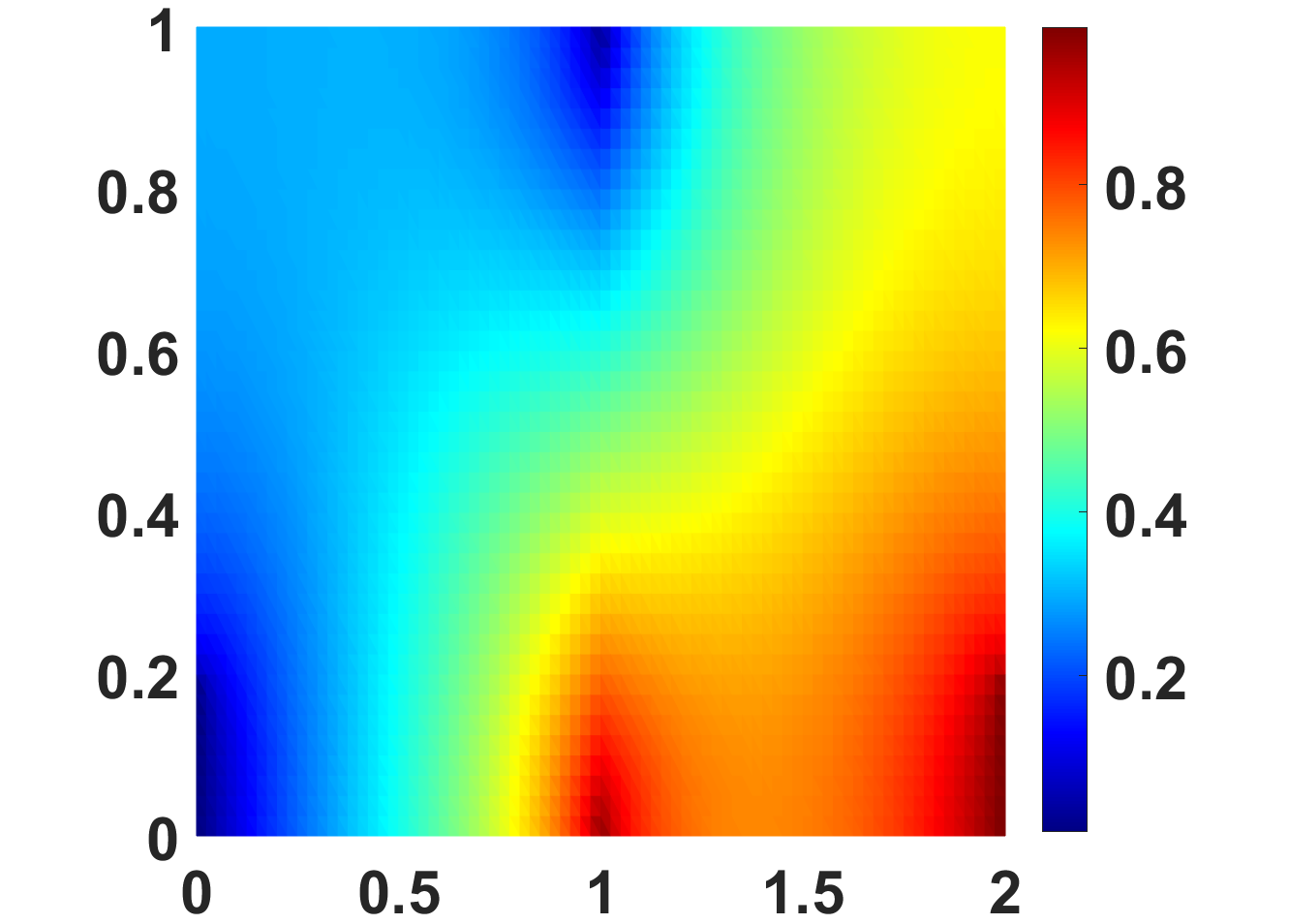}
  \label{fig:pressure_field_nonimmersed}
\end{minipage}%
\begin{minipage}{.47\textwidth}
  \centering
  \includegraphics[width=1\linewidth]{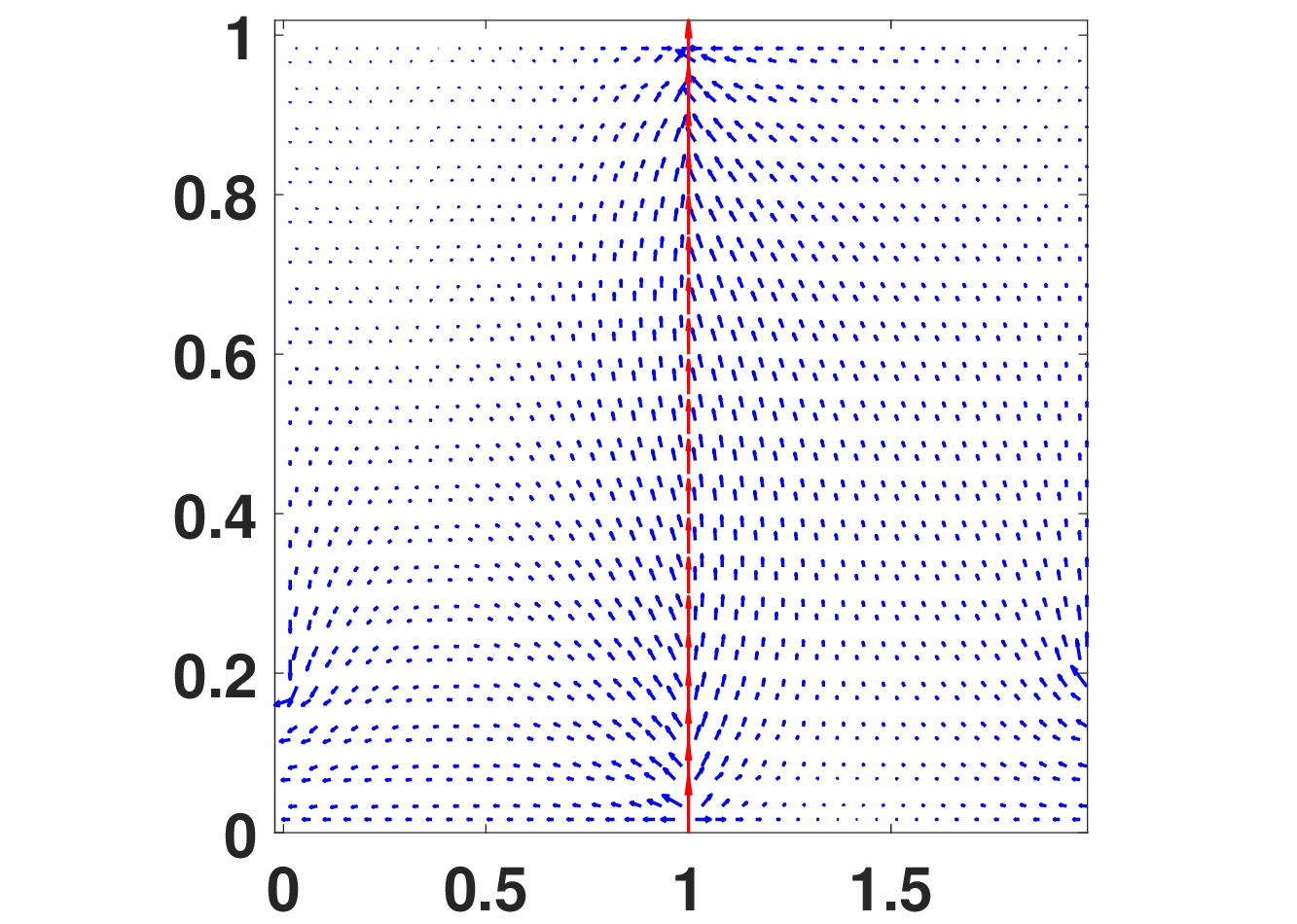}
  \label{fig:velocity_field_nonimmersed}
\end{minipage} \vspace{-0.2cm}
\caption{[Test case 1] Pressure field (left) and velocity field (right) at time $T=5$.}
\label{Pressure_Vel_fields_NonImmersed} \vspace{-0.3cm}
\end{figure}

For the first test case, the domain of calculation $\Omega=(0,2) \times (0,1)$ is divided into two equally sized subdomains by a fracture of width $\Toan{d = 0.001}$ parallel to the $y$-axis (see Figure~\ref{Non_Immersed}).  For the boundary conditions, we impose $p =1$ at the bottom and $p=0$ at the top of the fracture.  On the external boundaries of the subdomains, a no flow boundary condition is imposed except on the lower fifth (length 0.2) of both lateral sides where a Dirichlet condition is imposed: $p = 1$ on the right and $p = 0$ on the left. In Figure~\ref{Pressure_Vel_fields_NonImmersed} , the snapshots of the pressure and velocity fields at final time $T = 5$ are shown to illustrate the behavior of the solutions.
\begin{figure}[h!]
\centering
 \includegraphics[width=0.58\linewidth]{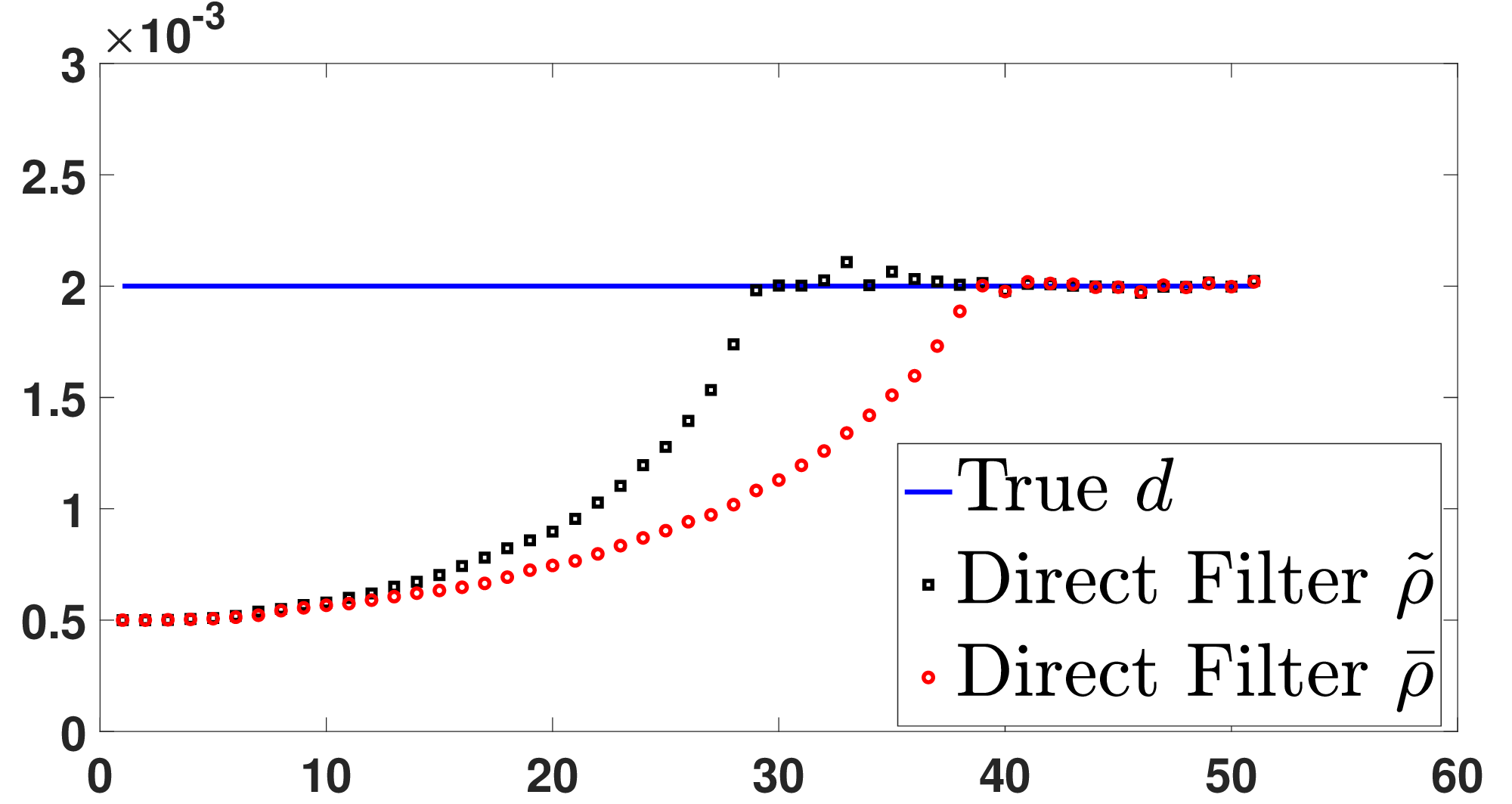}
  \label{fig:test2}
\caption{[Test case 1] The estimations for $d$.}
\label{DF_SingleFracture} 
\end{figure}

In this test case, we choose $M=80$ particles to empirically approximate the distribution of the unknown parameter and let $\xi_n$ in \eqref{ObserSetup_Egs} to be $\xi_n \sim \mathcal{N}(0, 500)$. We also consider two sets of noise $\bar{\epsilon} = \{\bar{\epsilon}_n\}$ and $\tilde{\epsilon} = \{\tilde{\epsilon}_n\}$ in \eqref{ParaSetup_Egs} where $\bar{\epsilon}_n \sim \mathcal{N}(0, 400)$ and $\tilde{\epsilon}_n \sim \mathcal{N}(0, 800)$. Finally, we denote $\bar{\rho} = \{1/\bar{\theta}_n\}$ and $\tilde{\rho} = \{1/\tilde{\theta}_n\}$ to be the sets of parameters obtained from the direct filter method corresponding to the noise $\bar{\epsilon}$ and $\tilde{\epsilon}$, respectively.

We present in Figure~\ref{DF_SingleFracture} the estimation for $d$ for each set of noise. We observe that both cases give accurate approximations for the true value $d$. Moreover, the direct filter $\Ph{\tilde{\rho}}$ tends to the true value of $\Toan{d}$ after $30$ steps, which is faster than the direct filter parameters $\Ph{\bar{\rho}}$ (around $40$ steps). This means that increasing the range of the noise, which is equivalent to increasing the exploration rate for the parameters, may accelerate the convergence speed of the algorithm.
\subsection{Test case 2: Two Parallel Fractures}
\begin{figure}[h!]
\centering
 \includegraphics[width=0.65\linewidth]{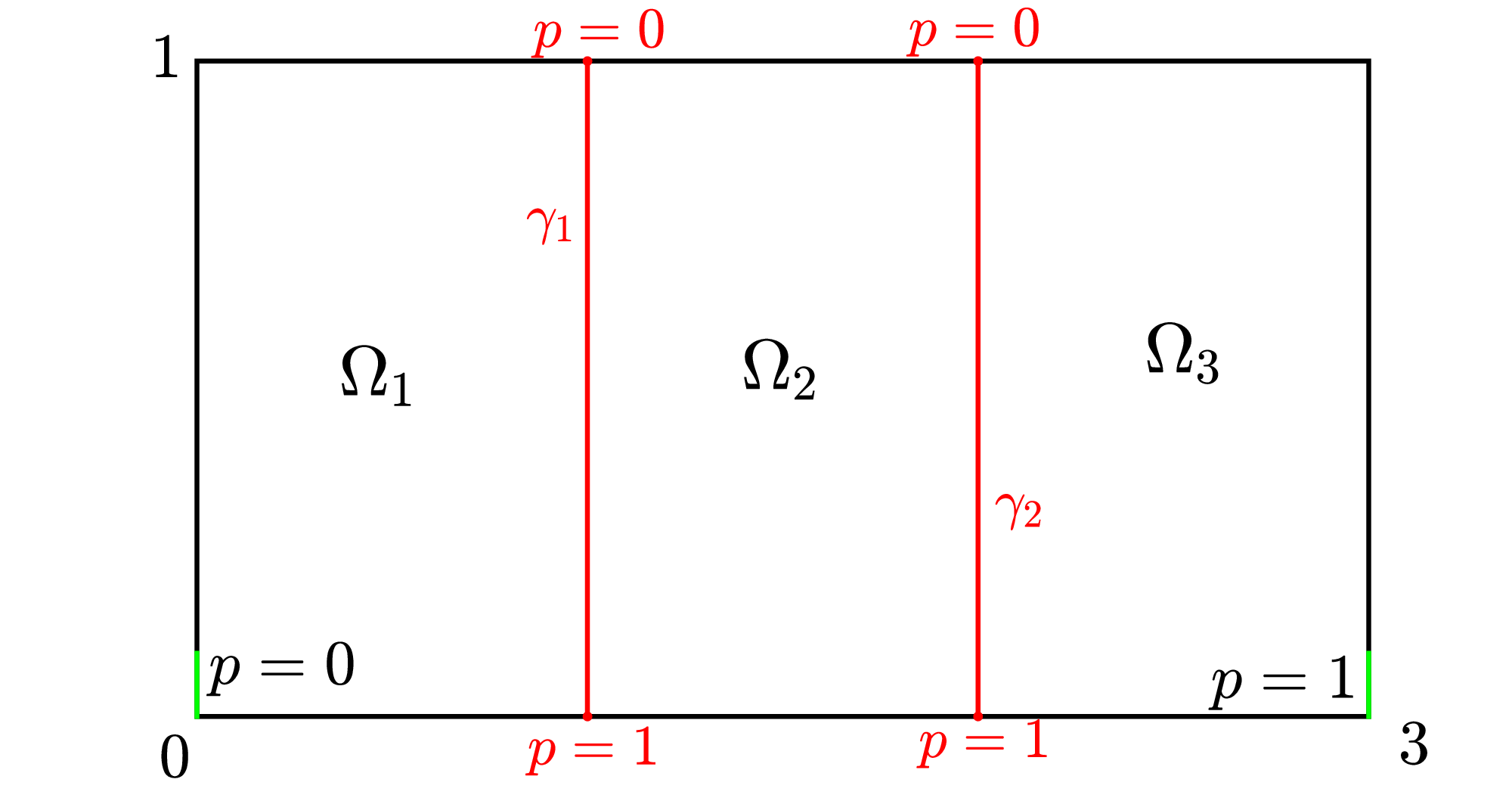}
  \label{fig:test2}
\caption{[Test case 2] Geometry and boundary conditions of the test case.}
\label{Non_Immersed_TwoFrac} \vspace{-0.1cm} 
\end{figure}%
In the second test case, we consider the domain of calculation being divided into three subdomains by two thin fractures $\gamma_1$ and $\gamma_2$, whose widths are $\Ph{d_1 =2.5e-03}$ and $\Ph{d_2=5e-03}$, respectively. The boundary conditions are the same as in Test case 1, with additional boundary conditions imposed on the second fracture $\gamma_2$ (see Figure~\ref{Non_Immersed_TwoFrac}). We note the the reduced fracture model \eqref{reduced_subdomain}-\eqref{reduced_fracture} and its discretization  \eqref{weak_discrete_reduced}-\eqref{initial_weak_discrete_reduced} can be extended straightforwardly to this case. 

\begin{figure}[h!]
\centering
\hspace{0.8cm}\begin{minipage}{.48\textwidth}\includegraphics[width=0.9\linewidth]{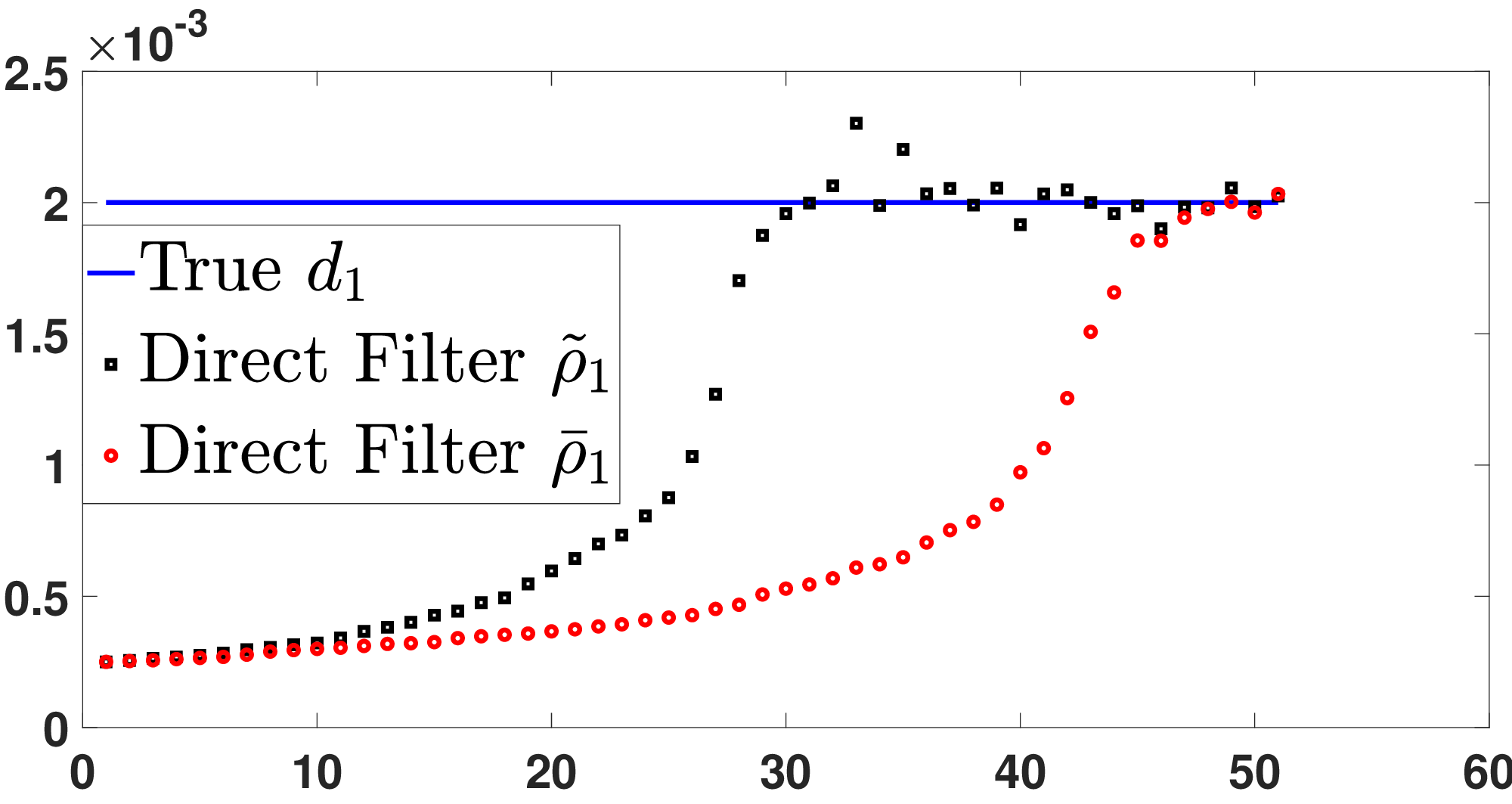}
  \label{fig:pressure_field_nonimmersed}
\end{minipage}
\begin{minipage}{.46\textwidth}
  \centering
\includegraphics[width=0.9\linewidth]{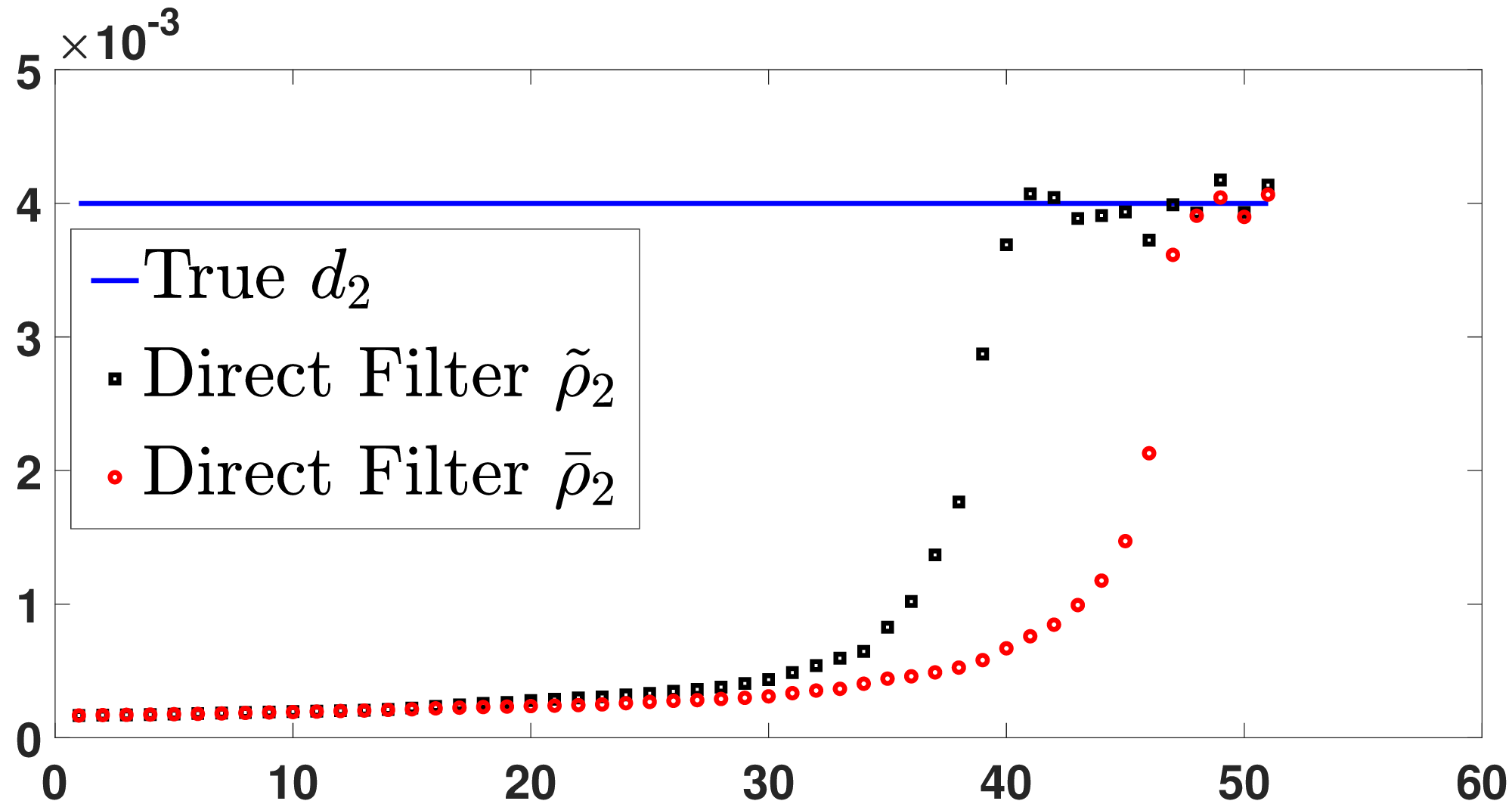}
  \label{fig:velocity_field_nonimmersed}
\end{minipage} \vspace{-0.4cm}
\caption{[Test case 2] The estimations for $d_1$ [left] and $d_2$ [right].}
\label{DFFracture_Parallel}
\end{figure}

We aim to approximate $\Toan{d_1}$ and $\Toan{d_2}$ using the direct filter algorithm. The parameter $\theta$ in \eqref{ParaSetup_Egs} is now a vector given by $\theta = (1/d_1, 1/d_2)$. We still choose the number of particles $M$ to be $80$ and $\xi_{n} \sim \mathcal{N}(0, 500)$. On the other hands, to accurately approximate both parameters, we need to choose bigger exploration rates. In particular, we consider two sets of noise $\bar{\epsilon}=\{\bar{\epsilon}_n\}_{n}=\{(\bar{\epsilon}_{1, n}, \bar{\epsilon}_{2, n})\}_n$ and  $\tilde{\epsilon} = \{\tilde{\epsilon}_n\}_n=\{(\tilde{\epsilon}_{1, n}, \tilde{\epsilon}_{2, n})\}_n$ where $\bar{\epsilon}_{k, n}$ and $\tilde{\epsilon}_{k, n}$ are the noises corresponding to the filtering particles for $d_k, \; k=1, 2$ with $\bar{\epsilon}_n\sim \mathcal{N}\left(0, \text{diag}(2000, 7000)\right)$ and $\tilde{\epsilon}_n\sim \mathcal{N}\left(0, \text{diag}(4000, 8000)\right)$. \Toan{Finally, we denote by $\bar{\rho} = \left(\bar{\rho}_1, \bar{\rho}_2\right)$ and $\tilde{\rho} =\left(\tilde{\rho}_1, \tilde{\rho}_2\right)$ the particle estimates corresponding to the noise $\bar{\epsilon}$ and $\tilde{\epsilon}$, respectively.}

We can observe from Figure~\ref{DFFracture_Parallel} that both cases provide nearly the same approximate values for the parameters of interested. Similar to Test case 1, bigger exploration rates lead to faster convergence speed. For example, it took 30 steps for $\tilde{\epsilon}$ to approach the true value of $d_1$, while nearly $47$ steps was required for $\bar{\epsilon}$ to reach the same value. \Ph{We also remark that increasing the exploration rates does not significantly raise the computational cost, so the method's efficiency remains preserved.} 
\subsection{Test case 3: Two intersecting fractures}
\begin{figure}[h!]
\centering
 \includegraphics[width=0.65\linewidth]{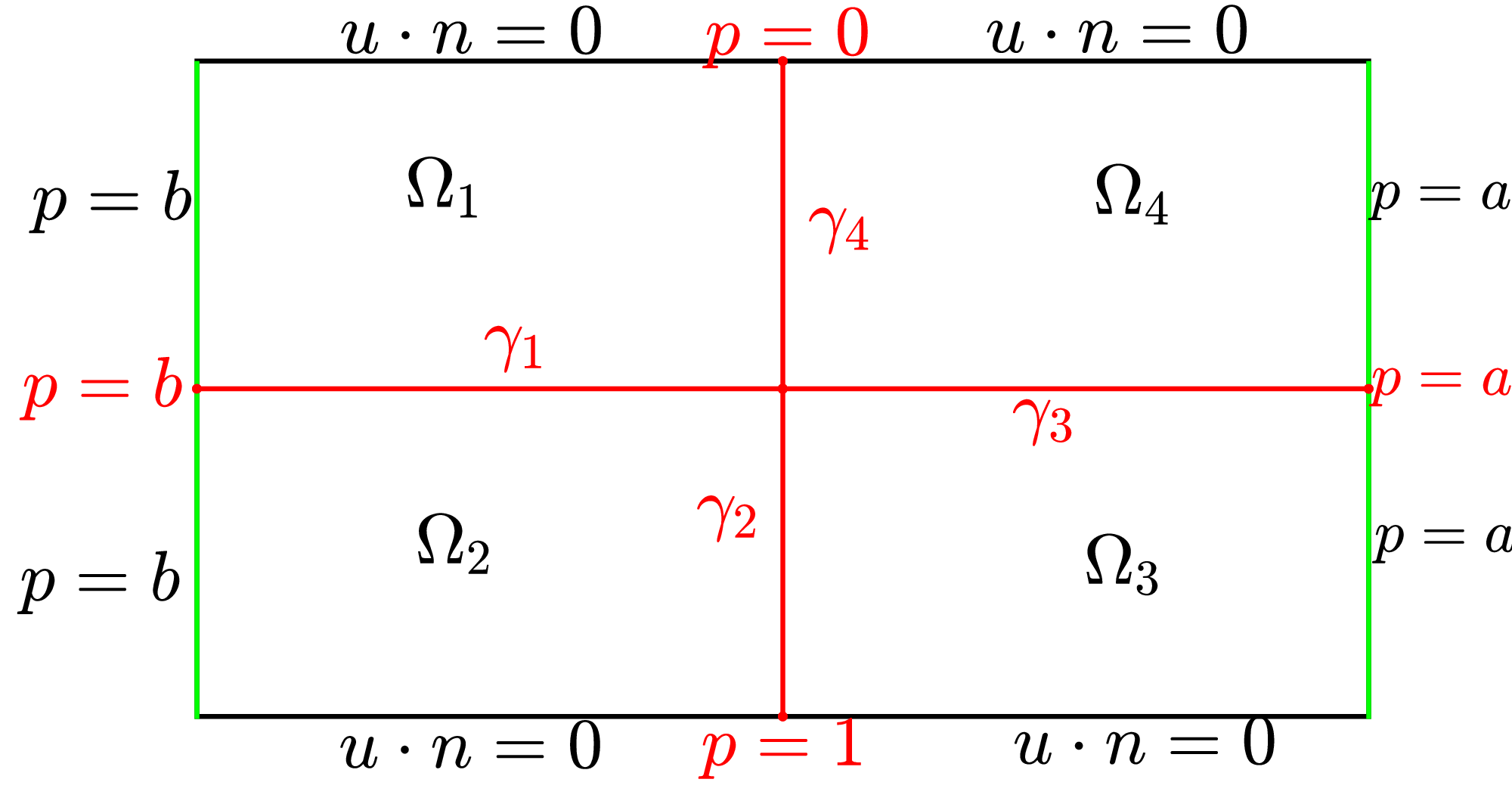}
  \label{fig:test2}
\caption{[Test case 3] Geometry and boundary conditions of the test case.}
\label{Non_Immersed_Intersecting} 
\end{figure}
We consider a test case adapted from~\cite{Alboin2002} where the domain of calculation is divided into four subdomains $\Omega_k, k=1,\hdots, 4$ by two intersecting fractures $\Gamma_1$ and $\Gamma_2$, whose widths are \Toan{$d_1=1e-03$ and $d_2 = 6e-04$}, respectively. To write the reduced model in this case, we further partition each fracture into two smaller parts: $\Gamma_1 = \gamma_1 \cup \gamma_3$ and $\Gamma_2 = \gamma_2 \cup \gamma_4$ (see Figure~\ref{Non_Immersed_Intersecting}). 
The boundary conditions are given in Figure~\ref{Non_Immersed_Intersecting}
.The reduced fracture model is then a coupled system between the four diffusion equations on $\Omega_k, k= 1, \hdots 4$ and the four tangential PDEs on $\gamma_k, k=1, \hdots, 4$. 

For each $k=1, \hdots, 4$, denote by $Z^n_{h, k} = \left(\bu^n_{h, k}, \; p^n_{h, k}\right)$ the solution on $\Omega_k$ and by $Z^n_{h, \gamma_k} = \left(\bu^n_{h, \gamma_k}, \; p^n_{h, \gamma_k}\right)$ the solution on $\gamma_k$ at time instant $n$. The solution $X_n$ at time instant $n$ is given by $X_n = \left(Z^n_{h, k}, \; Z^n_{h, \gamma_k}\right)_{k=1,\hdots, 4}$, and the noisy observational data can be written as
\begin{align*}
Y_{n+1} = \left(Z^n_{h, \gamma_k}\right)_{k=1,\hdots, 4} + \pmb{\xi}_{n+1} =\left(Z^n_{h, \gamma_k}\right)_{k=1,\hdots, 4} + \left(\xi^k_{n+1}\right)_{k=1,\hdots, 4} .
\end{align*}
The vector of targer parameters $\theta$ is set to be the same as Test case $2$, which is $\theta =(1/d_1, 1/d_2)$, and the noise $\epsilon_n$ in \eqref{ParaSetup_Egs} also consists of two components $\epsilon_n = \left(\epsilon_{1, n}, \epsilon_{2, n}\right)$. Unlike Test case 2, we need to imposed additional conditions at the intersection point of $\Gamma_1$ and $\Gamma_2$. Since the permeability in the fractures is considered to be \Toan{much larger} than the one in the subdomains, we can impose the following boundary conditions at the intersection point~(\cite{Amir2021, Alboin2002, Forma2014}):
\begin{align*}
\begin{array}{l}
p^n_{h, \gamma_1} = p^n_{h, \gamma_2} = p^n_{h, \gamma_3} =p^n_{h, \gamma_4}, \vspace{0.1cm} \\
 u^n_{h, \gamma_1}+u^n_{h, \gamma_2}+u^n_{h, \gamma_3}+u^n_{h, \gamma_4}=0.
\end{array}
\end{align*}
These conditions express the continuity of the \Toan{pressure} and the normal \Toan{fluxes} on the fractures across the intersection point.
\begin{figure}[h!]
\centering
\begin{minipage}{.5\textwidth}
  \centering
  \includegraphics[width=0.9\linewidth]{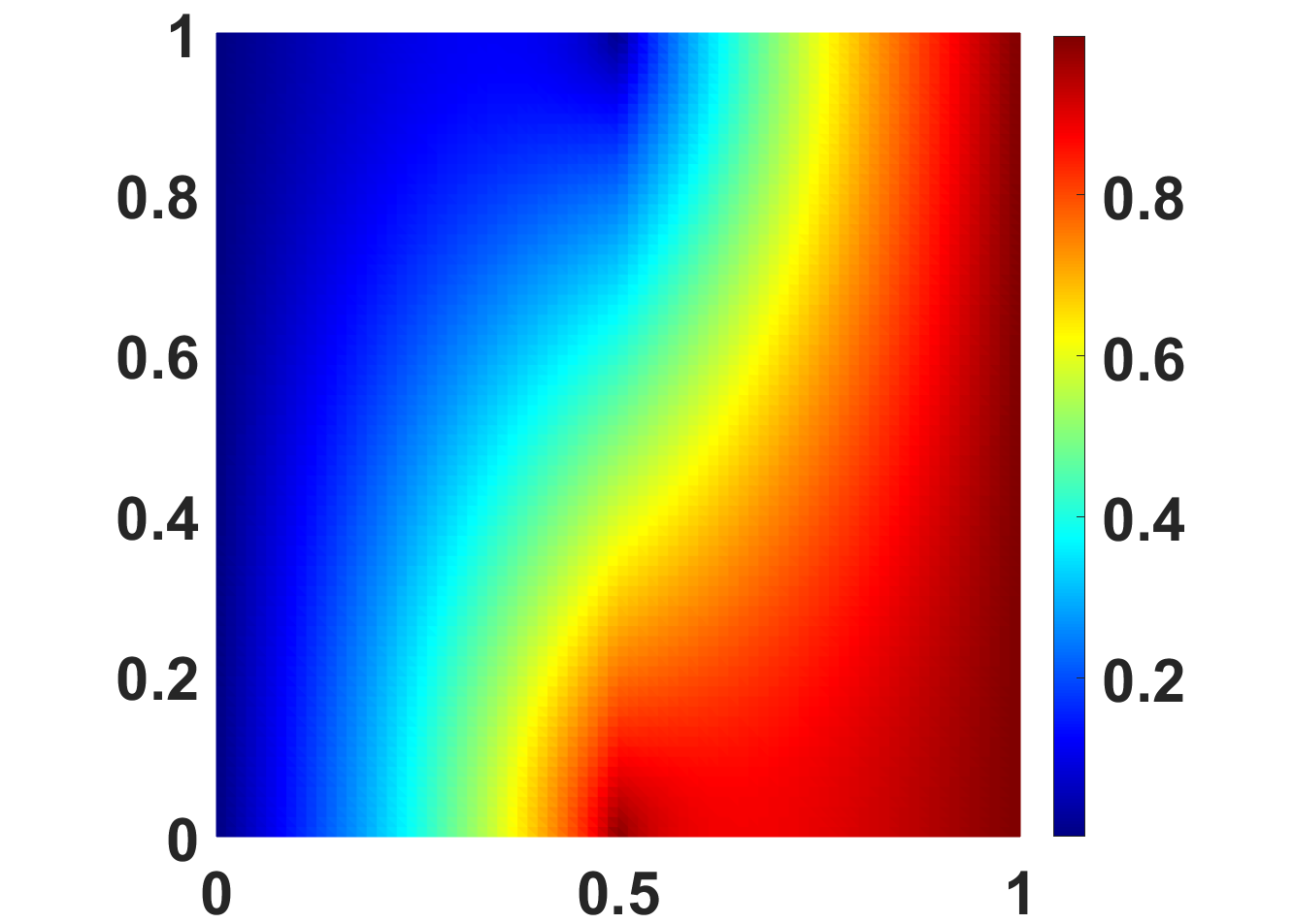}
  \label{fig:pressure_field_intersecting}
\end{minipage}%
\begin{minipage}{.5\textwidth}
  \centering
  \includegraphics[width=0.9\linewidth]{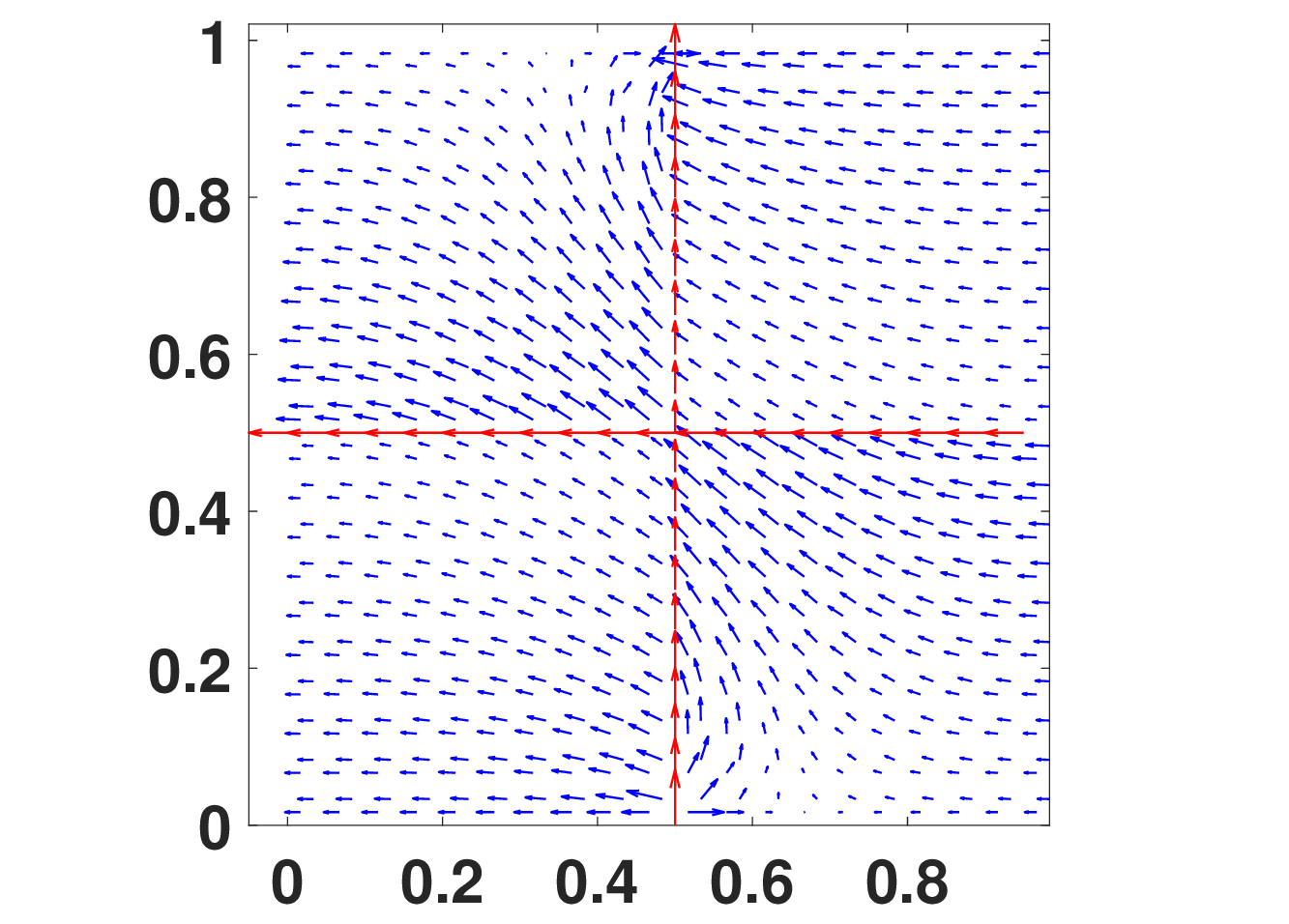}
  \label{fig:velocity_field_intersecting}
\end{minipage} 
\caption{[Test case 3a] Pressure field (left) and velocity field (right) for $a=1$ and $b=0$ at time $T=5$.}
\label{Pressure_Vel_fields_Intersecting_1st}
\end{figure}%
\begin{figure}[h!]
\centering
 \includegraphics[width=0.68\linewidth]{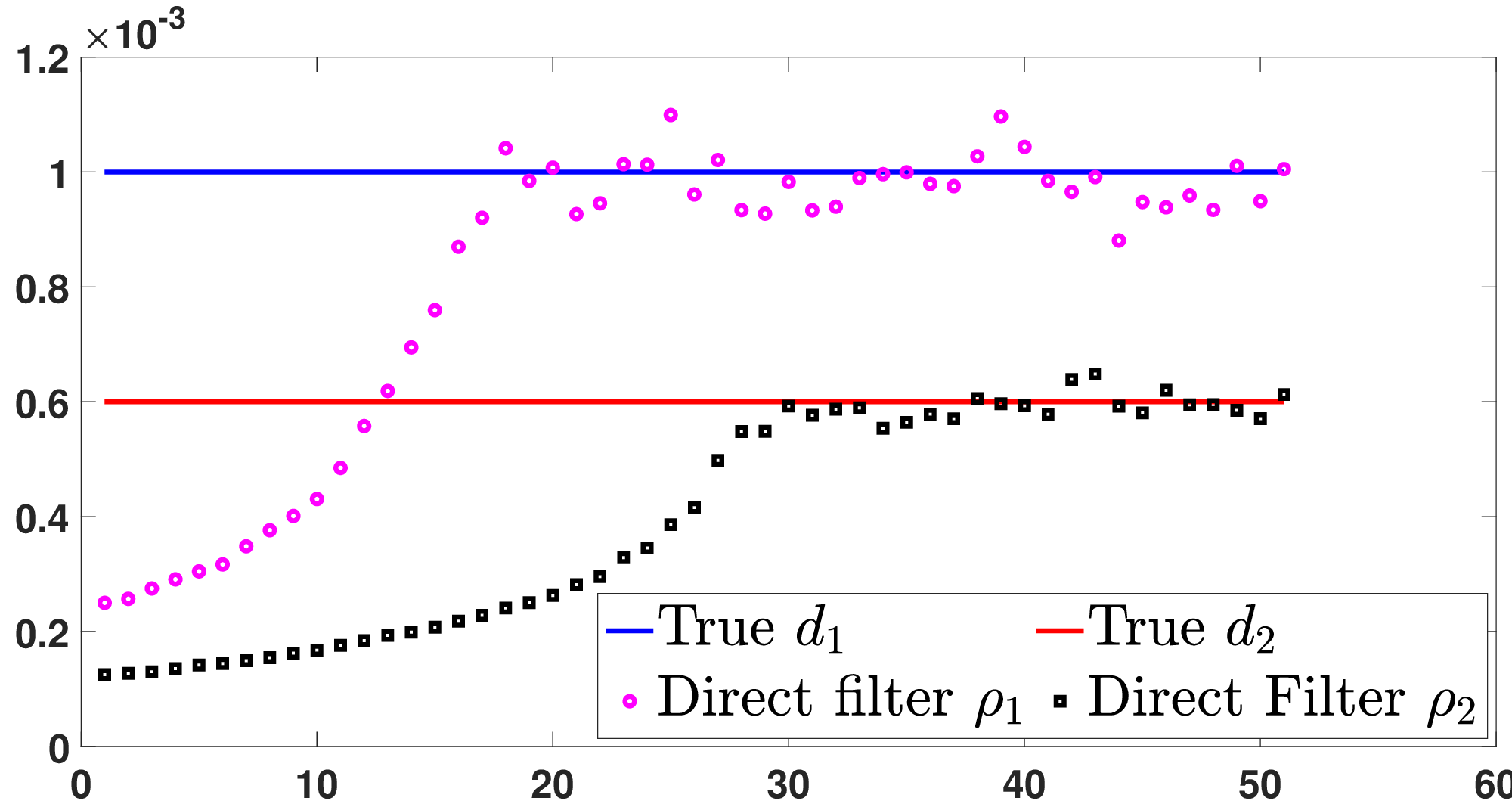}
  \label{fig:test2}
\caption{[Test case 3a] The estimations for $d_1$ and $d_2$ .}
\label{DFIntersecting_FirstBCs} 
\end{figure}

We consider two sets of boundary conditions for Test case 3. For the first set, we let $a=1$ and $b=0$, which leads to Test case 3a. We show the snapshot of the solution for Test case 3a in Figure~\ref{Pressure_Vel_fields_Intersecting_1st}. We can observe that the magnitude of the fluxes on the two fractures are relatively the same, and they \Ph{flow} in the direction driven by the boundary conditions. For Test case 3a, we fix the number of particles to be $M=120$ and  choose $\epsilon_{n} \sim \mathcal{N}(0, \text{diag}(8000, 10000))$ and present the direct filter estimations for $\rho_1 = 1/\theta_1$ and $\rho_2 = 1/\theta_2$ in Figure~\ref{DFIntersecting_FirstBCs}. With the chosen noise, we can see that the parameters converge to the true values relatively fast. \Ph{We observed that the rate of convergence of $\rho_2$ is slower than that of $\rho_1$. It can be explained that the flow in the vertical fracture is a combination between the vertical flow in the fracture and the flow from the right subdomain. This mixed information may cause the parameter $d_2$ to be sensitive, so that convergence rate is slower.}

\begin{figure}[h!]
\centering
\begin{minipage}{.5\textwidth}
  \centering
  \includegraphics[width=0.9\linewidth]{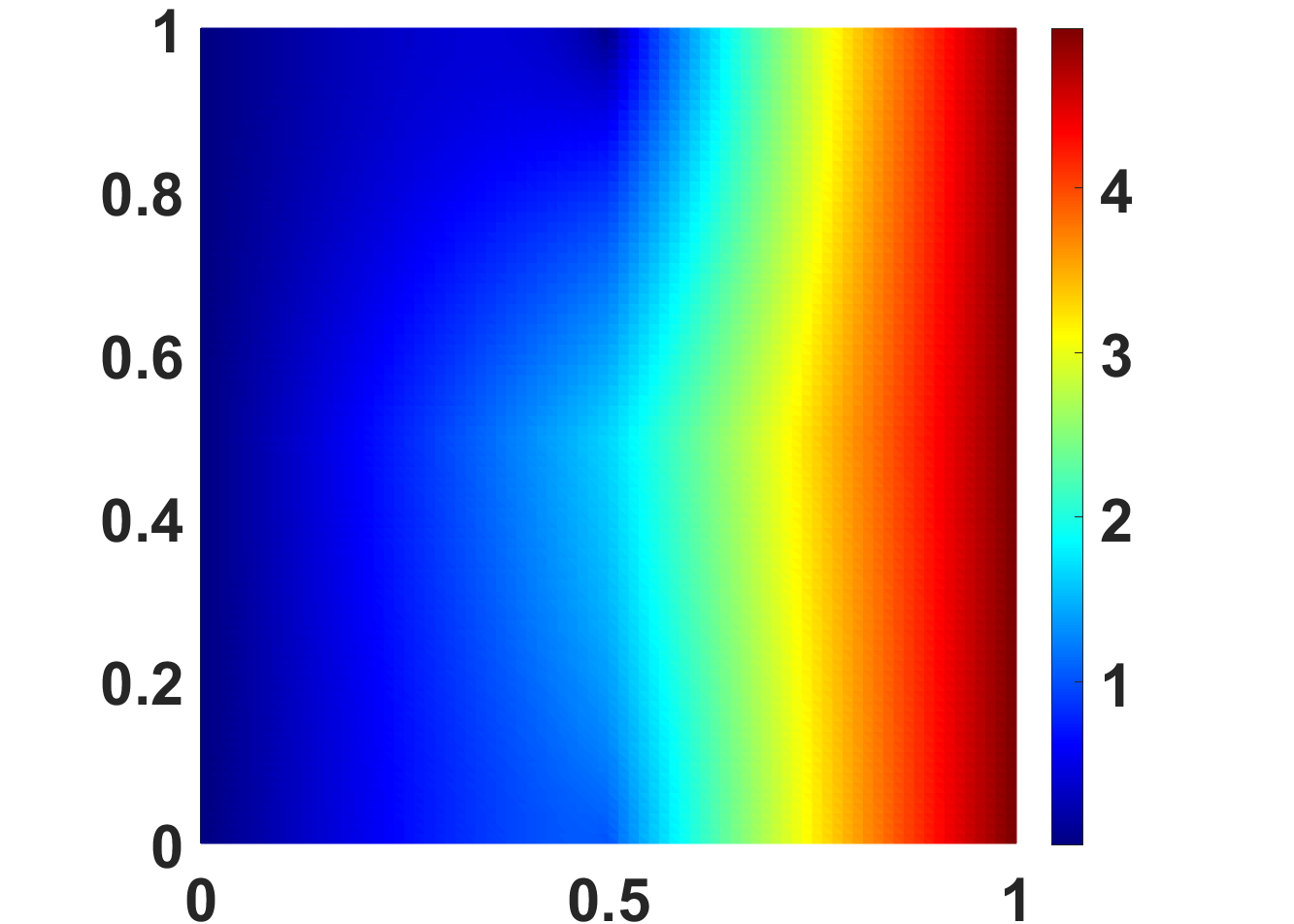}
  \label{fig:pressure_field_intersecting}
\end{minipage}%
\begin{minipage}{.5\textwidth}
  \centering
  \includegraphics[width=0.9\linewidth]{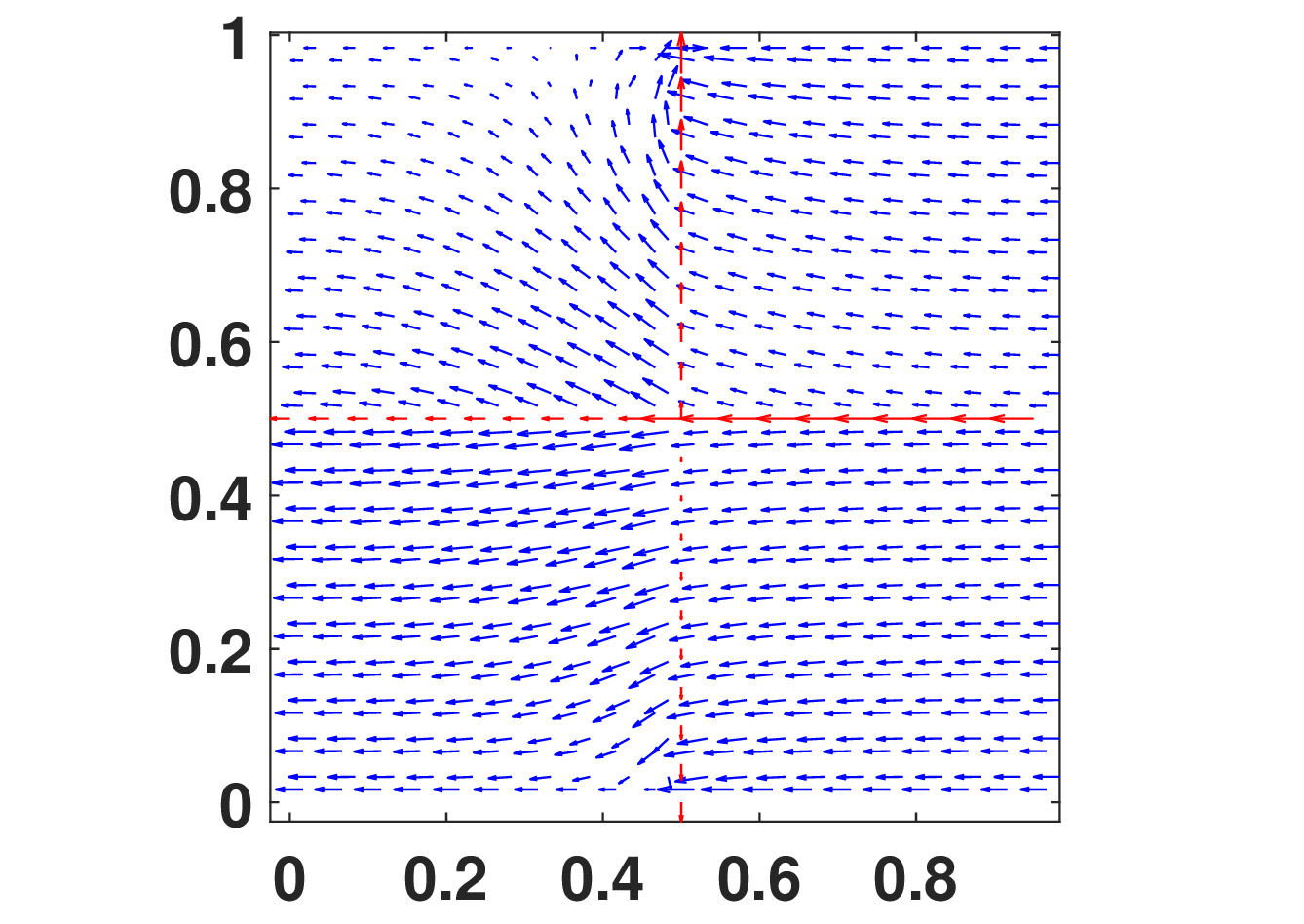}
  \label{fig:velocity_field_intersecting}
\end{minipage} 
\caption{[Test case 3b] Pressure field (left) and velocity field (right) for $a=5$ and $b=0$.}
\label{Pressure_Vel_fields_Intersecting_2nd} \vspace{-0.2cm}
\end{figure}

\begin{figure}[h!]
\centering
 \includegraphics[width=0.7\linewidth]{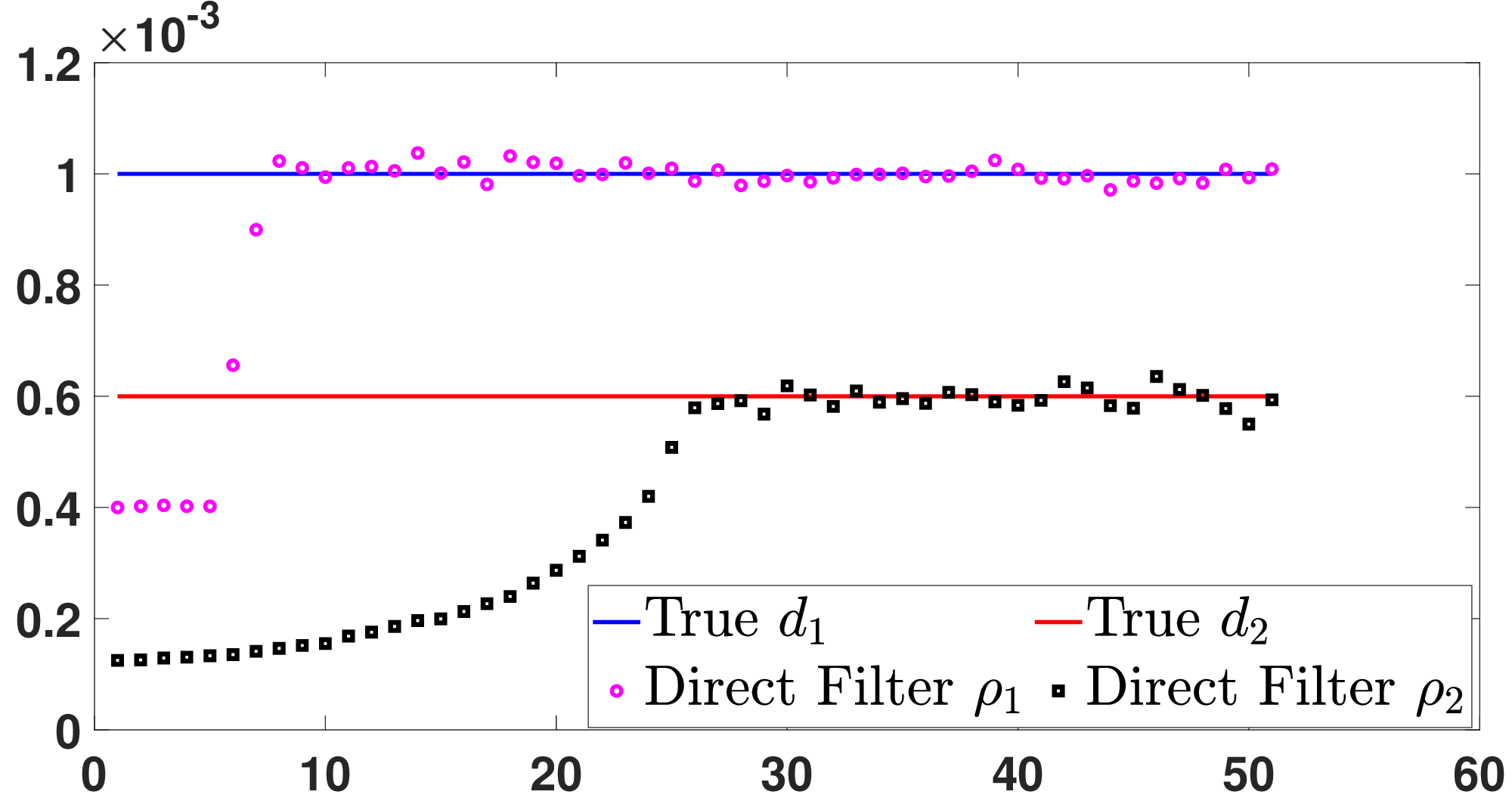}
  \label{fig:test2}
\caption{[Test case 3b] The estimations for $d_1$ and $d_2$.}
\label{DFIntersecting_SecondBCs} 
\end{figure}

Next, we consider a different set of boundary conditions where $a=5$ and $b=0$, which leads us to Test case 3b. The snapshot of the solution for this test case is shown in Figure~\ref{Pressure_Vel_fields_Intersecting_2nd}.We can see that the flux in the horizontal fracture dominates the vertical one, which drives the flux on the lower part of the vertical fracture to flow in the opposite direction compared to Test case 3a. For this test case, to accurately approximate the values of $d_1$ and $d_2$, we need to choose large values for the noise, which is $\epsilon_{n} \sim \mathcal{N}(0, \text{diag}(18000, 18000))$.  The \Ph{parameter estimates are presented} in Figure~\ref{DFIntersecting_SecondBCs}. \Toan{Again, we observe that the convergence rate of $\rho_2$ is slower than that of $\rho_1$ due to the sensitivity of the parameter $d_2$. However, similar to the previous test cases, the convergence speeds of the parameters to the exact values are relatively fast.}


\section{Conclusion}
In this work, we studied the inverse problem of estimating the widths of the fractures in a fractured porous medium based on observations of the fluid flow in the rock matrix. A reduced fracture model was introduced and was fully discretized to serve as the forward problem of the model inversion. The inverse problem was then solved by using an online parameter estimation technique by adopting the direct filter method developed in~\cite{Bao2019a}. We presented several numerical experiments to show that the parameters of interest can be recovered accurately by our method under various circumstances. For the future work, we aim to extend our work to a more general goal where a complete characterization of a fracture can be inferred based on the observational data in the porous medium.

\bibliographystyle{plain}

\bibliography{References}
\end{document}